\documentclass[jsl]{asl}
\usepackage[utf8x]{inputenc}
\usepackage{t1enc}

\usepackage[hidelinks]{hyperref}
\usepackage{proof}
\usepackage{amssymb}
\usepackage{txfonts}
\usepackage{url}
\usepackage{color}

\newcommand{\defeq}{\triangleq}
\newcommand{\imp}{\Rightarrow}
\newcommand{\myiff}{\Leftrightarrow}
\newcommand{\arrow}{\rightarrow}

\newcommand{\dimp}{\,\dot{\imp}\,}
\newcommand{\dneg}{\dot{\neg}}
\newcommand{\dbot}{\dot{\bot}}
\newcommand{\dand}{\,\dot{\wedge}\,}
\newcommand{\dor}{\,\dot{\vee}\,}
\newcommand{\dforall}{\dot{\forall}}
\newcommand{\dexists}{\dot{\exists}}

\newcommand \seqr[3]
  {\shortstack{$#2$ \\ \mbox{}\\
                   \mbox{}\hrulefill\mbox{}\\ \mbox{}\\ $#3$} \raisebox{2ex}{$\;\;\mbox{$#1$}$}}

\newcommand{\substid}{\mathit{id}}
\newcommand{\godel}[1]{\lceil #1 \rceil}
\newcommand{\truth}[3]{{#3} \vDash_{#2} #1}
\newcommand{\ptruth}[3]{{#3} \vDash^e_{#2} #1}
\newcommand{\pvDash}{\vDash^e}
\newcommand{\termtruth}[3]{[\![#1]\!]^{#2}_{#3}}

\newcommand{\DNS}{\mathsf{DNS}}
\newcommand{\DDNS}{\mathsf{DDNS}}
\newcommand{\EFQ}[1]{\mathit{EFQ}_{#1}}
\newcommand{\WFT}{\mathsf{WFT}}
\newcommand{\isprovable}[2]{#1 \vdash #2}
\newcommand{\consistent}[1]{#1 \not \vdash \dbot}
\newcommand{\consistentexpanded}[1]{(#1 \vdash \dbot) \imp \bot}

\newcommand{\metarulestyle}[1]{\mathit{#1}}
\newcommand{\dest}[3]{\metarulestyle{dest}\,#1\,\metarulestyle{as}\,#2\,\metarulestyle{in}~#3}

\newcommand{\destlinetwo}[3]{\!\!\!\begin{array}{l}\metarulestyle{dest}\,#1\,\\\metarulestyle{as}\,#2\,\\\metarulestyle{in}~#3\end{array}}

\newcommand{\destlinetopin}[3]{\!\!\!\begin{array}{l}\metarulestyle{dest}\,#1\,\metarulestyle{as}\,#2\,\metarulestyle{in}\\#3\end{array}}

\newcommand{\Ccaseline}[5]{\metarulestyle{case}~{#1}~\metarulestyle{of}~\!\!\!\begin{array}[t]{l}[{#2} \mapsto {#3}~\\|\,{#4} \mapsto {#5}]\end{array}}
\newcommand{\dcaseraw}{\objectrulestyle{\dot{case}}}
\newcommand{\dcase}[3]{\objectrulestyle{\dot{case}}\,#1\,#2\,#3}
\newcommand{\efq}{\metarulestyle{efq}}

\newcommand{\metafunctionstyle}[1]{\mathit{#1}}
\newcommand{\completeness}{\metafunctionstyle{completeness}}
\newcommand{\soundness}{\metafunctionstyle{soundness}}
\newcommand{\init}{\metafunctionstyle{init}}
\newcommand{\squeeze}{\metafunctionstyle{flush}}
\newcommand{\kontimp}[2]{\metafunctionstyle{kont}^{#1}_{#2}}
\newcommand{\kontimpprim}[2]{\metafunctionstyle{kont'}^{#1}_{#2}}

\newcommand{\kontorkont}[2]{\metafunctionstyle{kont}^{#1}_{#2}}
\newcommand{\Jbase}{\metafunctionstyle{J_{\metafunctionstyle{base}}}}
\newcommand{\Jcons}{\metafunctionstyle{J_{\metafunctionstyle{cons}}}}
\newcommand{\IO}{\metafunctionstyle{I_0}}
\newcommand{\IS}{\metafunctionstyle{I_S}}

\newcommand{\Iimp}{\metafunctionstyle{I_{\imp}}}
\newcommand{\Iforall}{\metafunctionstyle{I_{\forall}}}
\newcommand{\subsetjoin}{\metafunctionstyle{join}_\subset}
\newcommand{\join}{\metafunctionstyle{join}}
\newcommand{\hjoin}{\metafunctionstyle{hjoin}}
\newcommand{\merge}{\metafunctionstyle{share}}
\newcommand{\mergethree}{\metafunctionstyle{share_3}}

\newcommand{\inclctx}{\metafunctionstyle{incl}}
\newcommand{\inclcupctx}{\metafunctionstyle{incl'}}
\newcommand{\inj}{\metafunctionstyle{inj}}
\newcommand{\dinj}{\objectrulestyle{\dot{inj}}}
\newcommand{\class}{\metafunctionstyle{classic}}

\newcommand{\objectrulestyle}[1]{\mathit{#1}}
\newcommand{\ax}{\dot\objectrulestyle{ax}}
\newcommand{\appimp}{\dot\objectrulestyle{app^{\dimp}}}
\newcommand{\app}{\dot\objectrulestyle{app}}
\newcommand{\appf}{\dot\objectrulestyle{app^{\dforall}}}
\newcommand{\abs}{\dot\objectrulestyle{abs}}
\newcommand{\absimp}{\dot\objectrulestyle{abs^{\dimp}}}

\newcommand{\drinker}{\dot\objectrulestyle{drinker}}
\newcommand{\mydn}{\dot\objectrulestyle{dn}}
\newcommand{\dweak}{\dot\objectrulestyle{weak}}

\newcommand{\CASE}{\dot\objectrulestyle{CASE}}
\newcommand{\APPIMP}{\dot\objectrulestyle{APP^{\dimp}}}
\newcommand{\APPF}{\dot\objectrulestyle{APP^{\dforall}}}
\newcommand{\AXDRINKER}{\dot\objectrulestyle{AX}}
\newcommand{\AXINIT}{\dot\objectrulestyle{AX}^0_{\dneg A_0}}
\newcommand{\AXIMP}{\dot\objectrulestyle{AX}}
\newcommand{\DNABS}{\dot\objectrulestyle{DNABS}}
\newcommand{\LIFTBOT}{\dot\objectrulestyle{BOT}}
\newcommand{\DNLARGE}{\dot\objectrulestyle{DN}}
\newcommand{\PROJLEFT}[1]{\dot\objectrulestyle{PROJ}^{\dimp}_1}
\newcommand{\PROJRIGHT}[1]{\dot\objectrulestyle{PROJ}^{\dimp}_2}
\newcommand{\ANDPAIR}[1]{\dot\objectrulestyle{PAIR}}
\newcommand{\ANDPROJLEFT}[1]{\dot\objectrulestyle{PROJ}^{\land}_1}
\newcommand{\ANDPROJRIGHT}[1]{\dot\objectrulestyle{PROJ}^{\land}_2}
\newcommand{\ANDPROJBOTHSIDES}[2]{\dot\objectrulestyle{PROJ}^{\land}_{#1}}
\newcommand{\ORINJLEFT}[1]{\dot\objectrulestyle{INJ}_1}
\newcommand{\ORINJRIGHT}[1]{\dot\objectrulestyle{INJ}_2}
\newcommand{\ORINJBOTHSIDES}[2]{\dot\objectrulestyle{INJ}_i}

\newcommand{\classic}{\mathcal{C}\mathit{lassic}}
\newcommand{\exfalso}{\mathcal{E}\mathit{xfalso}}
\newcommand{\peirce}{\mathcal{P}\mathit{eirce}}

\newcommand{\Fun}{{\mathcal{F}\!un}}
\newcommand{\Term}{\mathcal{T}\!erm}
\renewcommand{\Form}{\mathcal{F}\!orm}
\newcommand{\Pred}{\mathcal{P}\!red}

\begin{document}

\title[The constructive content of Henkin's proof of Gödel's
  completeness theorem]
{An analysis of the constructive content of Henkin's proof of Gödel's
  completeness theorem}

\author{Hugo Herbelin}
\revauthor{Herbelin, Hugo}
\address{IRIF, CNRS, Inria, Universit\'e Paris Cit\'e}
\email{Hugo.Herbelin@inria.fr}

\author{Danko Ilik}
\revauthor{Ilik, Danko}
\address{CNES, France}
\email{dankoilik@gmail.com}

\begin{abstract}
G\"odel's completeness theorem for classical first-order logic is one of the
most basic theorems of logic. Central to any foundational course in
logic, it connects the notion of valid formula to the notion of
provable formula.

We survey a few standard formulations and proofs of the completeness theorem before focusing
on the formal description of a slight modification of Henkin's proof within
intuitionistic second-order arithmetic.

It is standard in the context of the completeness of intuitionistic
logic with respect to various semantics such as Kripke or Beth
semantics to follow the Curry-Howard correspondence and to interpret
the proofs of completeness as programs which turn proofs of validity
for these semantics into proofs of derivability.

We apply this approach to Henkin's proof to phrase it as a program
which transforms any proof of validity with respect to Tarski semantics
into a proof of derivability.

By doing so, we hope to shed an ``effective'' light on the relation
between Tarski semantics and syntax: proofs of validity are
syntactic objects that we can manipulate and compute with, like
ordinary syntax.
\end{abstract}

\maketitle

\tableofcontents
\vspace{-8mm}

\newpage
\section{Preliminaries}

\subsection{The completeness theorem}

The completeness theorem for classical first-order logic is one of the most
basic and traditional theorems of logic. Proved by Gödel in
1929~\cite{Godel29} as an answer to a question raised by Hilbert and
Ackermann in 1928~\cite{HilbertAckermann28}, it states that any of the
standard equivalent formal systems for defining provability in
first-order logic is complete enough
to include a derivation of every valid formula. A formula $A$
is {\em valid} when it is true under all interpretations of its
primitive symbols over any domain of quantification.

Let $\mathcal{L}$ be a signature for first-order logic, i.e. the data of a
set\footnote{We use here ``set'' in an informal way, not
  necessarily assuming the metalanguage to be specifically set
  theory.} $\Fun$ of function symbols, each of them coming with an arity, as well as of a set $\Pred$ of predicate symbols, each
of them also coming with an arity. We call constants the function
symbols of arity~0 and propositional atoms the predicate symbols of
arity~0. When studying the computational content
of Gödel's completeness in Section~\ref{sec:Henkin}, we shall restrict the
language to a countable one but the rest of this section does not
require restrictions on the cardinal of the language.

We let $f$ range over $\Fun$ and $P$ range over $\Pred$.  For
$f\in\Fun$ and $P\in\Pred$, we respectively write their arity $a_f$
and $a_P$. Let $x$ range over a countable set $\mathcal{X}$ of variables
and let $t$ range over the set $\Term$ of terms over $\mathcal{L}$ as
described by the following grammar:
$$
\begin{array}{lll}
t & ::= & x ~\mid~ f(t_1,...,t_{a_f})\\
\end{array}
$$
Let $A$ range over the set $\Form$ of formulae over $\mathcal{L}$
as described by the following grammar:
$$
\begin{array}{lll}
A & ::= & P(t_1,...,t_{P_f}) ~\mid~ \dbot ~\mid~ A \dimp A ~\mid~ \dforall x\,A\\
\end{array}
$$
Note that in classical first-order logic, the language of {\em negative} connectives and quantifiers made of
$\dimp$, $\dforall$ and $\dbot$ is enough to express all other connectives and quantifiers.
The dot over the notations is to distinguish the connectives and quantifier of the logic we are
talking about (object logic) from the connectives and quantifiers of the ambient logic in which
the completeness theorem is formulated (meta-logic, see below). We take $\dbot$ as a
primitive connective and this allows to express consistency of the object logic as the non-provability of $\dbot$. Negation can then be defined as
$\dneg A \defeq A \dimp \dbot$. Also, in $\dforall x\,A$, we say that
$x$ is a binding variable which binds all occurrences of $x$
in $A$ (if any). If the occurrence of a variable is not in the scope of
a $\dforall$ with same name, it is called free. If a formula has no
free variables, we say it is {\em closed}.

Let us write $\Gamma$ for finite contexts of hypotheses, as defined by
the following grammar:
$$
\begin{array}{lll}
\Gamma & ::= & \epsilon ~\mid~ \Gamma, A\\
\end{array}
$$ In particular, $\epsilon$ denotes the empty context, which we might
also not write at all, as e.g. in $\vdash A$ standing for $\epsilon
\vdash A$.

We assume having chosen a formal system for provability in classical
first-order logic, e.g. one of the axiomatic systems
 given in Frege~\cite{Frege1879} or in Hilbert
and Ackermann~\cite{HilbertAckermann28}, or one of the systems such as
Gentzen-Jaśkowski's natural deduction~\cite{Jaskowski34,Gentzen35} or
Gentzen's sequent calculus~\cite{Gentzen35}, etc., and we write
$\isprovable{\Gamma}{A}$ for the statement that $A$ is
provable under the finite context of hypotheses $\Gamma$. If $\mathcal{M}$
 is a model for classical logic
and $\sigma$ an interpretation of the variables from $\mathcal{X}$
 in the model, we write $\truth{A}{\sigma}{\mathcal{M}}$ for the statement
expressing that $A$ is true in the model $\mathcal{M}$ (to be defined in Section~\ref{sec:models}). \emph{Validity} of~$A$ under
assumptions $\Gamma$, written $\Gamma \vDash A$ is defined to be
$\forall \mathcal{M}\;\forall \sigma\;(\truth{\Gamma}{\sigma}{\mathcal{M}}\imp \truth{A}{\sigma}{\mathcal{M}})$
where $\truth{\Gamma}{\sigma}{\mathcal{M}}$ is the conjunction of all
$\truth{B}{\sigma}{\mathcal{M}}$ for every $B$ in $\Gamma$,
i.e. $\bigwedge_{B\in\Gamma} \truth{B}{\sigma}{\mathcal{M}}$. Note that $\imp$,
$\forall$, $\wedge$ and, later on, below, $\lor$, $\exists$, $\bot$, as
well as derived $\neg$, represent the connectives and quantifiers of
the metalanguage.

We say that \emph{$\Gamma$ is inconsistent} if $\Gamma \vdash \dbot$ and
\emph{consistent} if $\consistent{\Gamma}$, i.e. if
$\consistentexpanded{\Gamma}$, i.e., if a contradiction in the object
language is reflected as a contradiction in the meta-logic. We say that
\emph{$\Gamma$ has a model} if there exist $\mathcal{M}$ and $\sigma$ such that
$\truth{\Gamma}{\sigma}{\mathcal{M}}$. The completeness theorem, actually a
weak form of the completeness theorem as discussed in the next
section, is commonly stated under one of the following classically but
not intuitionistically equivalent forms:
$$
\begin{array}{ll}
C1. & \vDash A ~\imp~\isprovable{}{A}\\
C2. & \mbox{$\Gamma$ is consistent} ~\imp~ \mbox{$\Gamma$ has a model}\\
C3. & \mbox{$\Gamma,\neg A$ has a model} ~\lor~ \isprovable{\Gamma}{A}\\
\end{array}
$$

\subsection{Weak and strong completeness}

In a strong form, referred to as strong completeness\footnote{We follow
here a terminology
dubbed by Henkin in his 1947 dissertation, according
to~\cite{Henkin96}.
  However, in
  the context of intuitionistic logic, some
  authors use the weak and strong adjectives with different
  meanings. For instance, in Kreisel~\cite{Kreisel58b,Kreisel62}, the statement
  $(\vDash A) \imp (\vdash A)$ is called strong completeness while
  weak completeness is the statement $(\vDash A) \imp \neg\neg(\vdash
  A)$.
\label{page:Okada}
  In the context of semantic cut-elimination, e.g. in
  Okada~\cite{Okada02}, $(\vDash A) \imp (\vdash A)$ is only a weak form
  of completeness whose strong form is the statement $(\vDash A) \imp
  (\vdash_{\mathrm{cut-free}} A)$, for a notion of cut-free proof similar
  to the notion of cut-free proof in Gentzen's sequent calculus or to normal proofs in Prawitz' analysis of normalisation for natural deduction.},
completeness states that any formula valid under some possibly
infinite theory is provable under a finite subset of this theory. This
is the most standard formulation of completeness in textbooks, and, as
such, it is a key component of the compactness theorem. Also proved by
Gödel~\cite{Godel30}, the compactness theorem states that it is enough
for a theory to have a model that any finite subset of the theory has
a model. In contrast, completeness with respect to finite theories as
stated above is referred to as weak completeness. Let $\mathcal{T}$ be a
set of formulae and let $\mathcal{T} \vdash A$ mean the existence of a
finite sequence $\Gamma$ of formulae in $\mathcal{T}$ such that $\Gamma
\vdash A$. Let $\mathcal{M} \vDash_\sigma \mathcal{T}$ be
$\forall B\in \mathcal{T}\,\truth{B}{\sigma}{\mathcal{M}}$ and let
the definitions of $\mathcal{T}$ is
consistent and of $\mathcal{T}$ has model be extended accordingly.  The strong
formulations of the three views at weak completeness above are now the
following:

$$
\begin{array}{ll}
S1. & \mathcal{T} \vDash A ~\imp~\isprovable{\mathcal{T}}{A}\\
S2. & \mbox{$\mathcal{T}$ is consistent} ~\imp~ \mbox{$\mathcal{T}$ has a model}\\
S3. & \mbox{$\mathcal{T} \cup \{\neg A\}$ has a model} ~\lor~ \isprovable{\mathcal{T}}{A}\\
\end{array}
$$

We shall consider the formalisation and computational content of
strong completeness. Weak completeness will then come as a
special case.

\subsection{The standard existing proofs of completeness}

Let us list a few traditional proofs from the classic
literature\footnote{We cite the most common proofs in the classic
  pre-1960 literature. Recent developments include e.g. Joyal's
  categorical presentation of a completeness theorem. We can also cite
  Berger's~\cite[Sec. 1.4.3]{SchwichtenbergWainer06} or
  Krivtsov~\cite{Krivtsov10} construction in intuitionistic logic of a classical model from a
  Beth model for classical provability. These two latter proofs are variants
  of the Beth-Hintikka-Kanger-Schütte style of proofs, the
  first one relying on the axiom of dependent choice and the second
  on the (weaker) Fan theorem.}.

\begin{itemize}

\item Gödel's original proof~\cite{Godel29} considers formulae in
  prenex form and works by induction on the number of quantifiers for
  reducing the completeness of first-order predicate logic
  completeness to the completeness of propositional logic.

\item Henkin's proof~\cite{Henkin49a} is related to statement S2: from
  the assumption that $\mathcal{T}$ is consistent, a syntactic model over
  the terms is built as a maximal consistent extension of $\mathcal{T}$
  obtained by ordering the set of formulae and extending
  $\mathcal{T}$ with those formulae that preserve
  consistency, following the ordering.

\item In the 1950's, a new kind of proof independently credited to
  Beth~\cite{Beth55}, Hintikka~\cite{Hintikka55,Hintikka55b},
  Kanger~\cite{Kanger57} and
  Sch\"utte~\cite{Schutte56} was given.
  The underlying idea is to build an
  infinite normal derivation, typically in sequent calculus. Rules are
  applied in a fair way, such that all possible combinations of rules
  are considered. If the derivation happens to be finite, a proof is
  obtained. Otherwise, by weak Kőnig's lemma, there is an infinite
  branch and this infinite branch gives rise to a countermodel. The
  intuition underlying this proof is then best represented by
  statement S3.

\item In the 1950's also, Rasiowa and
  Sikorski~\cite{RasiowaSikorski50} gave a variant of Henkin's proof
  relying on the existence of an ultrafilter for the Lindenbaum algebra
  of classes of logically equivalent formulae, identifying validity
  with having value 1 in all interpretations of a formula within
  the two-value Boolean algebra $\{0,1\}$. This is close to
  Henkin's proof in the sense that Henkin's proof implicitly
  builds an ultrafilter of the Lindenbaum algebra of formulae.

\end{itemize}

Our main contribution in this paper is the analysis in Section~\ref{sec:Henkin} of the
computational content of Henkin's proof.

\subsection{Models and truth}
\label{sec:models}
The interpretation of terms in a model $\mathcal{M}$ is given by a domain
$\mathcal{D}$ and by an interpretation $\mathcal{F}$ of
the symbols in $\Fun$ such that $\mathcal{F}(f) \in \mathcal{D}^{a_f} \arrow \mathcal{D}$,
where $\mathcal{D}^{a_f} \arrow \mathcal{D}$
denotes the set of functions of
arity $a_f$ over $\mathcal{D}$. Then, given an assignment
$\sigma\in \mathcal{X} \arrow \mathcal{D}$ of the variables to
arbitrary values of the domain, the interpretation of terms in $\mathcal{D}$ is given by:
$$
\begin{array}{lll}
\termtruth{x}{\sigma}{\mathcal{M}} & \defeq & \sigma(x)\\
\termtruth{f(t_1,\ldots,t_{a_f})}{\sigma}{\mathcal{M}} & \defeq & \mathcal{F}(f) (\termtruth{t_1}{\sigma}{\mathcal{M}}, \ldots, \termtruth{t_{a_f}}{\sigma}{\mathcal{M}})\\
\end{array}
$$

To interpret formulae, two common approaches are used in the
literature.

\begin{itemize}
\item {\em Tarski semantics (predicates as predicates).} This is the approach followed e.g. in
  the Handbook of Mathematical Logic~\cite{HandbookBarwise}, the
  Handbook of Proof Theory~\cite{HandbookBuss} or in the original
  proof of G\"odel~\cite{Godel29}. This approach interprets formulae
  of the object language propositionally, i.e. as formulae of the
  metalanguage. In this case, the interpretation depends on whether
  the metalanguage is classical or not. For instance, in a classical
  metalanguage,
  the theory $$\classic ~\defeq~ \{\dneg \dneg A \dimp A~|~A\in
  {\Form}\}$$ would be true in all models. On the other hand, in an
  intuitionistic metalanguage, a formula such as, say, $\dneg\dneg X \dimp
  X$ 
  could not be proved true in all models\footnote{For instance, if
    $\mathit{coh_{ML}}$ is the formula expressing the consistency of
    the metalanguage represented as an object language in the
    metalanguage itself, then a model $\mathcal{M}$ binding atom $X$ to
    the metalanguage formula $\mathit{coh_{ML}} \vee \neg
    \mathit{coh_{ML}}$ would intuitionistically satisfy $\mathcal{M}
    \vDash \neg\neg X$ but not $\mathcal{M} \vDash X$.}. In a
  strongly anti-classical intuitionistic metalanguage refuting
  double-negation elimination, it could even be proved that there are
  models\footnote{For instance, in second-order intuitionistic
    arithmetic extended with Church Thesis (CT), excluded-middle on
    undecidable formulae is provably contradictory and the same model
    interpreting $X$ as $\mathit{coh_{ML}} \vee \neg
    \mathit{coh_{ML}}$ invalidates $\dneg\dneg X \dimp X$.} which
  refute $\dneg\dneg X \dimp X$.

  The possible presence of models provably anti-classical is
  not a problem per se for proving completeness as completeness is only about exhibiting
  one particular model and it is possible to ensure that $\dneg\dneg X \dimp X$
  holds in this particular model. However, whether the metalanguage is classical or not has an
  impact on the soundness property, i.e. on the statement that the
  provability of $A$ implies the validity of $A$. Indeed, there is little hope
  to prove the soundness of double-negation elimination if the quantification over
  models include non-classical models.
  Therefore, for the definition of validity to be both sound and
  complete for classical provability with respect to Tarski semantics, independently
  of whether the metalanguage is intuitionistic or classical, we
  would need to define classical validity using an explicit restriction to classical
  models:
  $$\mathcal{T} \vDash A ~\defeq ~ \forall \mathcal{M}\;\forall
  \sigma\;(\truth{\classic}{\sigma}{\mathcal{M}} \imp
  \truth{\mathcal{T}}{\sigma}{\mathcal{M}}\imp \truth{A}{\sigma}{\mathcal{M}})$$

\item {\em Bivalent semantics (predicates as binary functions).} Another approach is
  to assign to formulae a truth value in the two-valued set $\{0,1\}$
  and to define $\truth{A}{\sigma}{\mathcal{M}}$ as
  $\mathit{truth}_{\mathcal{M}}(A,\sigma) = 1$ for the
  corresponding $\mathit{truth}$ function. This is the
  approach followed e.g. in Rasiowa-Sikorski's proof, or also
  e.g. in~\cite{Church56,Simpson99}, among others\footnote{For
    instance, the definition of validity used in
    Henkin~\cite{Henkin49a}, though not fully formal, also intends a
    two-valued semantics.}. In particular, relying on a two-valued
  truth makes the theory $\classic$ automatically true.

  Depending on the metalanguage, a function from $\Form$ to $\{0,1\}$
  can itself be represented either as a functional relation,
  i.e. as a relation $\mathit{istrue}$ on $\Form \times \{0,1\}$ such
  that for all $A$, there is a unique $b$ such that
  $\mathit{istrue}(A,b)$ holds (this is the representation used e.g. for the
  completeness proof in~\cite{Simpson99}), or, primitively as a
  function if ever the metalanguage provides such primitive notion of
  function (as is typically the case in intuitionistic logics, e.g.
  Heyting Arithmetic in finite types~\cite{Troelstra73}, or
  Martin-Löf's type theory~\cite{MartinLof75,CoqPau89}).

  Reverse mathematics of the subsystems of classical second-order
  arithmetic have shown that building a model from a proof of
  consistency requires the full strength of $\Sigma^0_1$-separability,
  or equivalently, of Weak Kőnig's Lemma~\cite{Simpson99}. This
  implies that the corresponding $\mathit{truth}$ function is in general not
  recursive~\cite{Kleene52}. Expecting $\mathit{truth}$ to be definable
  primitively as a computable function in an intuitionistic logic is
  thus hopeless. As for representing truth by a functional
  relation $\mathit{istrue}$, the expected property
  $\mathit{istrue}(A,0) \lor \mathit{istrue}(A,1)$ could only be
  proven by requiring some amount of classical reasoning.

It is known how to compute with classical logic in second-order
arithmetic~\cite{Griffin90,Parigot92,Krivine94} and we could study the
computational content of a formalisation of the completeness proof
which uses this definition of truth. The extra need for classical
reasoning in this approach looks however like a useless complication,
so we shall concentrate on the predicates-as-predicates approach.

\end{itemize}

So, to summarise, we will not expect truth to be two-valued and will
require explicitly as a counterpart that models are classical,
leading to the following refined definitions\footnote{For the record, note that,
  in the presence of only negative connectives, an equivalent way to
  define $\vDash A$ so that it means the same in an intuitionistic and
  classical setting is to
  replace the definition of $\truth{P(t_1,...,t_{a_P})}{\sigma}{\mathcal{M}}$ by
  $$\truth{P(t_1,...,t_{a_P})}{\sigma}{\mathcal{M}} \defeq \neg\neg (\termtruth{t_1}{\sigma}{\mathcal{M}}, \ldots, \termtruth{t_{a_P}}{\sigma}{\mathcal{M}}) \in \mathcal{M}(P)$$
  or even, saving a negation as in Krivine~\cite{Krivine90b}, by
  $$\truth{P(t_1,...,t_{a_P})}{\sigma}{\mathcal{M}} \defeq \neg (\termtruth{t_1}{\sigma}{\mathcal{M}}, \ldots, \termtruth{t_{a_P}}{\sigma}{\mathcal{M}}) \in \mathcal{M}(P)$$
Indeed, in these cases, the definition of truth becomes a purely
negative formula for which intuitionistic and classical provability
coincide.} of validity and existence of a model:
$$
\begin{array}{lll}
\mathcal{T} \vDash A & \defeq & \forall \mathcal{M}\, \forall \sigma\, (\truth{\classic}{\sigma}{\mathcal{M}} \imp \truth{\mathcal{T}}{\sigma}{\mathcal{M}} \imp \truth{A}{\sigma}{\mathcal{M}})\\
\mbox{$\mathcal{T}$ has a model} & \defeq & \exists \mathcal{M}\, \exists \sigma\, (\truth{\classic}{\sigma}{\mathcal{M}} \wedge \truth{\mathcal{T}}{\sigma}{\mathcal{M}})\\
\end{array}
$$

Two auxiliary choices of presentation of Tarski semantics can be made\footnote{These auxiliary choices would have been relevant as well if we had chosen to represent truth as a map to $\{0,1\}$.}.

\begin{itemize}
\item {\em Recursively-defined truth.} The approach followed e.g.
  in the Handbook of Mathematical
  Logic~\cite{HandbookBarwise} or the Handbook of Proof
  Theory~\cite{HandbookBuss} is to have the model interpret only the
  predicate symbols and to have the truth of formulae defined recursively. This is obtained by giving an interpretation
$\mathcal{P}$ where any symbol $P\in\Pred$ is mapped to a set $\mathcal{P}(P) \subset \mathcal{D}^{a_P}$. Then, the
truth of a formula with respect to some assignment $\sigma$ of the
free variables is given recursively by:
$$
\begin{array}{lll}
\truth{P(t_1,\ldots,t_{a_f})}{\sigma}{\mathcal{M}} & \defeq & (\termtruth{t_1}{\sigma}{\mathcal{M}}, \ldots, \termtruth{t_{a_P}}{\sigma}{\mathcal{M}}) \in \mathcal{P}(P)\\
\truth{\dbot}{\sigma}{\mathcal{M}} & \defeq & \bot\\
\truth{A \dimp B}{\sigma}{\mathcal{M}} & \defeq & \truth{A}{\sigma}{\mathcal{M}} \imp \truth{B}{\sigma}{\mathcal{M}}\\
\truth{\dforall x\, A}{\sigma}{\mathcal{M}} & \defeq & \forall v\in\mathcal{D}\,\truth{A}{\sigma\cup[x\leftarrow v]}{\mathcal{M}}\\
\end{array}
$$
\label{sec:extensional-truth}

\item {\em Axiomatically-defined truth.} A common alternative approach is to
  define truth as a subset
  $\mathcal{S}$ of closed formulae in the language of terms
  extended with the constants of $\mathcal{D}$, such that:
  $\dbot$ is not in $\mathcal{S}$;
  $A \dimp B$ is in
  $\mathcal{S}$ iff $B$ is whenever $A$ is; $\dforall x\,A$
  is in $\mathcal{S}$ iff $A[x\leftarrow v]$ is for all values
  $v\in \mathcal{D}$; $A[x\leftarrow f(v_1,...,v_n)]$ is in
  $\mathcal{S}$ iff $A[x\leftarrow v]$ is in $\mathcal{S}$
  whenever $\mathcal{F}(f)(v_1,...,v_n) = v$ for some value
  $v \in \mathcal{S}$. This approach is adopted e.g. in Krivine~\cite{Krivine96b}.

\end{itemize}

We will retain the first approach which conveniently exempts us from
defining the set of formulae enriched with constants from $\mathcal{D}$.
So, shortly, a model $\mathcal{M}$ will be a triple $(\mathcal{D},
\mathcal{F}, \mathcal{P})$ where $\mathcal{F}$ maps any symbol $f\in\Fun$ to a
function $\mathcal{F}(f) \in \mathcal{D}^{a_f} \arrow \mathcal{D}$ and $\mathcal{P}$
maps any symbol $P\in\Pred$ to a set $\mathcal{P}(P) \subset \mathcal{D}^{a_P}$.

\subsection{Regarding the metalanguage as a formal system}

Let $M$ be the metalanguage in which completeness is stated and $O$ be
the object language used to represent provability in first-order
logic. In $M$, a proof of the validity of a formula $A$ is essentially
a proof of the universal closure of $A$, seen as a formula of $M$,
with the closure made over the domain of quantification of quantifiers, over the free
predicate symbols, over the free function symbols and over the free
variables of $A$. Otherwise
said, adopting a constructive view at proofs of the metalanguage,
we can think of the weak completeness theorem in form C1 as a process to transform a
proof of the universal closure of $A$ expressed in $M$ into a proof of
$A$ expressed in the proof object language $O$ (and conversely, the soundness
theorem can be seen as stating an embedding of $O$ into $M$).
Similarly, a proof of the validity of a formula $A$ with respect to an
infinite theory $\mathcal{T}$ is a proof in $M$ of the universal closure
of~$(\forall B \in \mathcal{T} \termtruth{B}{}{}) \imp \termtruth{A}{}{}$ where $\termtruth{C}{}{}$
 is the replication of $C$ as a formula of $M$ and,
 computationally speaking, statement S1 is a process to
turn such a proof in $M$ (which has to use only a finite subset of
$\mathcal{T}$ in $M$, since $M$, seen itself as a formal system, supports
only finite proofs) into a proof in~$O$.

The key point is however that this transformation of a proof in $M$
into a proof in $O$ is done in $M$ itself, and, within $M$ itself,
the only way to extract information out of a proof of validity is by
instantiating the free symbols of the interpretation of $A$ in $M$ by
actual function and predicate symbols of $M$, i.e. by producing what
at the end is a model, i.e. a domain, functions
and predicates actually definable in $M$.

\subsection{The computational content of
completeness proofs for intuitionistic logic}
\label{sec:semantic-normalisation}

It is known that composing the soundness and completeness theorems for
propositional or predicate logic gives a cut-elimination theorem, as
soon as completeness is formulated in such a way that it produces a
normal proof\footnote{Using e.g. Beth-Hintikka-Kanger-Sch\"utte's proof.}.
Now, if the proofs of soundness and completeness are formalised in a metalanguage
equipped with a normalisation procedure, e.g. in a $\lambda$-calculus-based proofs-as-programs presentation of
second-order arithmetic~\cite[Ch. 9]{Krivine90,Krivine93}, one gets an effective
cut-elimination theorem, namely an effective procedure which turns any
non-necessarily-normal proof of $\Gamma \vdash A$ into a normal proof of
$\Gamma \vdash A$.

In the context of intuitionistic provability, this has been explored
abundantly under the name of {\em semantic normalisation}, or {\em
  normalisation by evaluation}. Initially based on ideas from Berger
and Schwichtenberg~\cite{BergerSchwichtenberg91} in the context of
simply-typed $\lambda$-calculus, it was studied for the realisability
semantics of second-order implicative propositional
logic\footnote{I.e., Girard-Reynolds System~$F$.} by Altenkirch,
Hofmann and Streicher~\cite{AltenkirchHofmannStreicher96}, for the
realisability semantics of implicative propositional logic in Hilbert
style\footnote{I.e., equivalently, simply-typed combinatory logic.} by
Coquand and Dybjer~\cite{CoquandDybjer97}, for the Kripke semantics of
implicative propositional logic in natural deduction
style\footnote{I.e., equivalently, simply-typed $\lambda$-calculus.} by
C.~Coquand~\cite{Coquand02}, for Heyting algebras by Hermant and
Lipton~\cite{Hermant05,HermantLipton10}, etc. It has also been applied
to phase semantics of linear logic by Okada~\cite{Okada02}. It also
connects to a normalisation technique in computer science called
Typed-Directed Partial Evaluation (TDPE)~\cite{Danvy96}.

Let us recall how this approach works in the case of minimal implicative
propositional logic (Figure~\ref{fig:object-minimal-logic})
using soundness and completeness with respect to Kripke models~\cite{Coquand02}. Let
$\mathcal{K}$ range over Kripke models $(\mathcal{W},\leq,\Vdash_X)$ where
$\leq$ is a preorder on $\mathcal{W}$ and $\Vdash_X$ a monotonic predicate
over $\mathcal{W}$ for each propositional atom $X$. Let $w$
range over $\mathcal{W}$, i.e.  worlds in the corresponding Kripke models. Let
us write $ w \Vdash_\mathcal{K} A$ (resp. $w \Vdash_\mathcal{K} \Gamma$) for
truth of $A$ (resp. for the conjunction of the truth of all formulae
in $\Gamma$) at world $w$ in the Kripke model $\mathcal{K}$. In
particular, $w \Vdash_\mathcal{K} A$ is extended from atoms to all
formulae by defining $w \Vdash_\mathcal{K} A \dimp B \defeq\forall w'
(w' \geq w \imp w' \Vdash_\mathcal{K} A \imp w' \Vdash_\mathcal{K} B)$.
Let us write
$\Gamma \vDash_I A$ for the validity of $A$ relative to $\Gamma $ at all
worlds of all Kripke models, i.e. for the formula $\forall \mathcal{K}
\forall w\; (w \Vdash_\mathcal{K} \Gamma \imp w \Vdash_\mathcal{K} A)$.

The metalanguage being here a $\lambda$-calculus, we shall write
its proofs as mathematical functions.
We write $x \mapsto t$ for the proof of an implication
as well as for the proof of a universal quantification, possibly also writing $(x:A)
\mapsto t$ to make explicit that $x$ is the name of a proof of $A$.
We shall represent modus ponens and instantiation of universal
quantification by function application, written $t\,u$. We shall use the
notation $()$ for the canonical proof of a nullary conjunction and the
notation $(t,u)$ for the proof of a binary conjunction, seen as a product
type and obtained by taking the pair of the proofs of the components
of the conjunction. To give a name $f$ to the proof of a statement of
the form $\forall x_1,...,x_n\,(A \imp B)$ we shall use the notation
$f^{x_1,...,x_n}(a:A) : B$ followed by clauses of the form
$f^{x_1,...,x_n}(a) \defeq t$ (for readability, we may also write some
of the $x_i$ as subscripts rather than superscripts of $f$).

For instance, the proof that Kripke forcing
is monotone, i.e. that $\forall w w'\,(w' \geq w \land w \Vdash A \imp
w' \Vdash A)$, can be written as the following function $\Uparrow_A$, recursive in the structure of $A$,
taking as arguments two worlds $w$ and $w'$:
$$\begin{array}{lcrlcllcl}
  \Uparrow_A^{w,w'} & : & \!\!\!\! & w' \geq w & \land & w \Vdash A & \!\!\!\! & \imp & w' \Vdash A\\
  \Uparrow_X^{w,w'} & & (\!\!\!\!& h & , & m &\!\!\!\!) &\defeq& p_X(h,m)\\
  \Uparrow_{A\dimp B}^{w,w'} & & (\!\!\!\!& h & , & m &\!\!\!\!)& \defeq& w'' \mapsto (h':w''\geq w') \mapsto m\,w''\,(\mathit{trans}(h,h'))\\
\end{array}$$
where $p_X$ is the proof of monotonicity of $\Vdash_X$ and
$\mathit{trans}$ is the proof of transitivity of $\geq$, both coming
with the definition of Kripke models, while, in the definition, $h$ is
a proof of $w' \geq w$ and $m$ a proof of $w \Vdash A$.

Similarly, the extension of $\Uparrow$ to a proof that forcing of
contexts is monotone can be written as follows, where we reuse the
notation $\Uparrow$, now with a context as index, to denote a proof of
$\forall w w'\,(w' \geq w \land w \Vdash \Gamma \imp w' \Vdash
\Gamma)$:
$$\begin{array}{lcrlcllcl}
  \Uparrow_\Gamma^{w,w'} & : & \!\!\!\!& w' \geq w & \land & w \Vdash \Gamma &\!\!\!\!& \imp & w' \Vdash \Gamma\\
  \Uparrow_\epsilon^{w,w'} & & (\!\!\!\!& h & , & () &\!\!\!\!)&\defeq& ()\\
  \Uparrow_{\Gamma,A}^{w,w'} & & (\!\!\!\!& h & , & (\sigma,m) &\!\!\!\!)&\defeq& (\Uparrow_{\Gamma}^{w,w'}(h,\sigma),\Uparrow_A^{w,w'}(h,m))\\
\end{array}$$

Let us write $\Gamma \vdash_I A$ for intuitionistic provability.
Let us consider the canonical proof $\soundness^{\Gamma}_A$ of
$(\Gamma \vdash_I A) \imp (\Gamma \Vdash_I A) $ proved by
induction on the derivation of $\Gamma \vdash_I
A$. We write the proof as a recursive function, recursively on the structure of formulae:
$$
\!\!\begin{array}{llcl}
\soundness^\Gamma_A ~~:\!\!\! & \Gamma \vdash_I A & \!\!\imp\!\! & \Gamma \vDash_I A\\
\soundness^\Gamma_A    & \ax_i & \defeq & \mathcal{K} \mapsto w \mapsto \sigma \mapsto \sigma(i)\\
\soundness^\Gamma_{A\imp B} \!\!& \abs(p) & \defeq & \!\!\!\begin{array}[t]{l}\mathcal{K} \mapsto w \mapsto \sigma \mapsto w' \mapsto (h:w\leq w') \mapsto m \mapsto \\ \soundness^{\Gamma,A}_B \,p\,\mathcal{K}\, w'\, (\Uparrow_{\Gamma}^{w,w'} (h,\sigma),m)\end{array}\\
\soundness^\Gamma_B & \app(p,q) & \defeq & \!\!\!\begin{array}[t]{l}\mathcal{K} \mapsto w \mapsto \sigma \mapsto \\ (\soundness^{\Gamma}_{A\imp B}\,p\, \mathcal{K}\, w\, \sigma)~w~\mathit{refl}~(\soundness^{\Gamma}_{A}\,q\, \mathcal{K}\, w\, \sigma)\end{array}\\
\end{array}
$$ where $u$ is a proof of $\Gamma \vdash_I A$ in the last line,
$\app$, $\abs$, $\ax_i$ are the name of inference rules defining
object-level implicational propositional logic in a natural deduction style (see
Figure~\ref{fig:object-minimal-logic}); $\sigma(i)$ is the
$(i+1)^{\mbox{\scriptsize{th}}}$ component of $\sigma$ starting from
the right, and $\mathit{refl}$ is the proof of reflexivity of $\geq$
coming with the definition of Kripke models.

\begin{figure}
\begin{center}
{\em Primitive rules}
\end{center}

$$\seqr{\ax_i}{|\Gamma'|=i}{\Gamma, A, \Gamma' \vdash_I A}
\qquad
\seqr{\app}{\Gamma \vdash_I A \dimp B \qquad \Gamma \vdash_I A}{\Gamma \vdash_I B}
\qquad
\seqr{\abs}{\Gamma, A \vdash_I B}{\Gamma \vdash_I A \dimp B}
$$
\smallskip

\begin{center}
{\em Admissible rule}
\end{center}

$$\qquad\seqr{\dweak}{\Gamma\subset \Gamma' \qquad \Gamma \vdash_I A}{\Gamma' \vdash_I A}$$
\caption{Inference rules characterising minimal implicational logic}
\label{fig:object-minimal-logic}
\end{figure}

Let us also consider the following somehow canonical proof of cut-free
completeness, $\completeness : (\Gamma \vDash_I A) \imp (\Gamma
\vdash_I^{\mathit{cf}} A)$. It is based on the universal model of context
$\mathcal{K}_0$ defined by taking for $\mathcal{W}$ the set of contexts $\Gamma$ ordered by
inclusion and $\Gamma \vdash_I^{\mathit{cf}}
X$ for the forcing $\Vdash_X$ of atom $X$ at world $\Gamma$. Now, the proof
proceeds by showing the two directions of $\Gamma
\Vdash_{\mathcal{K}_0} A ~{\Leftrightarrow}~ \Gamma
\vdash_I^{\mathit{cf}} A$ by mutual induction on $A$. It is common to write
$\downarrow$ for the left-to-right direction (called \emph{reify}, or \emph{quote})
and $\uparrow$ for the right-to-left direction (called \emph{reflect}, or
\emph{eval}):
$$
\begin{array}{llcl}
{\downarrow^\Gamma_A~:} & \Gamma \Vdash_{\mathcal{K}_0} A  & {\imp} & \Gamma \vdash_I^{\mathit{cf}} A\\
\downarrow^\Gamma_{P} & m & \defeq & m\\
\downarrow^\Gamma_{A \dimp B} & m & \defeq & \abs\,(\downarrow^{\Gamma,A}_B (m \;(\Gamma,A) \;\inj^A_\Gamma \; (\uparrow^{\Gamma,A}_A\, \ax_0)))\\
\\
{\uparrow^\Gamma_A~:} & \Gamma \vdash_I^{\mathit{cf}} A & {\imp} & \Gamma \Vdash_{\mathcal{K}_0} A\\
\uparrow^\Gamma_{P} & p & \defeq & p \\
\uparrow^\Gamma_{A \dimp B} & p & \defeq & \Gamma' \mapsto f \mapsto m \mapsto \uparrow^{\Gamma'}_B(\app(\dweak(f,p),\downarrow^{\Gamma'}_A m))\\
\\
\init^\Gamma_{\Gamma'} & & : & \Gamma \Vdash_{\mathcal{K}_0} \Gamma'\\
\init^\Gamma_\epsilon  & & \defeq & ()\\
\init^\Gamma_{\Gamma',A} & & \defeq & (\init^{\Gamma}_{\Gamma'},\uparrow^\Gamma_A(\ax_{|\Gamma|-|\Gamma',A|}))\\
\\
\completeness^\Gamma_A~ : & \Gamma \vDash_I A & \imp & \Gamma \vdash_I^{\mathit{cf}} A\\
\completeness^\Gamma_A  & m & \defeq & \downarrow^\Gamma_A~ (m~\mathcal{K}_0~\Gamma~\init^\Gamma_\Gamma)
\end{array}
$$ where $|\Gamma|$ is the length of
$\Gamma$, $\dweak$ is the admissible rule of weakening in object-level
implicational propositional logic and $\inj^A_\Gamma$ is a
proof of $\Gamma \subset \Gamma,A$.

In particular, by placing ourselves in a metametalanguage, such that the
metalanguage is seen as a proofs-as-programs-style natural deduction
object language, i.e. as a typed $\lambda$-calculus, one would be able to show that

\begin{itemize}
\item for every given proof of $\Gamma \vdash_I A$, soundness produces,
  by normalisation of the metalanguage\footnote{Typically proved by embedding
  in another language assumed to be consistent.}, a proof of $\Gamma \vDash_I A$
  whose structure follows the one of the proof of $\Gamma \vdash_I A$;
\item for every proof of validity taken in canonical form (i.e. as a
  closed $\beta$-normal $\eta$-long $\lambda$-term of type $\Gamma
  \vDash A$ in the metalanguage), the resulting proof of $\Gamma
  \vdash_I^{\mathit{cf}} A$ obtained by completeness is, by
  normalisation in the metalanguage, the same $\lambda$-term with
  the abstractions and applications over $\mathcal{K}$, $w$ and proofs of
  $w \leq w'$ removed.
\end{itemize}

On the other side, if our proofs-as-programs-based metalanguage is
able to state properties of its proofs (as is the case for instance of
Martin-L\"of's style type theories~\cite{MartinLof75}), it can be
shown within the metalanguage itself that the composition of
completeness and soundness produces normal forms. This is what
C.~Coquand did by showing that the above proofs of soundness and
completeness, seen as typed programs, satisfy the following properties:
$$\begin{array}{l}
\forall p : (\Gamma \vdash_I A)\;p \sim \soundness^\Gamma_A\,p\,\mathcal{K}_0\,\Gamma\, init^{\Gamma}_{\Gamma} 
\forall p :(\Gamma \vdash_I A)\,\forall m :(\Gamma \vDash_I A)\; (p \sim m\,\mathcal{K}_0\,\Gamma\, init^{\Gamma}_{\Gamma} ~\imp~ p =_{\beta\eta} \completeness^{\Gamma}_A\,m)\\
\end{array}
$$ where
$\sim$ is an appropriate ``Tait computability'' relation
between object proofs and semantic proofs expressing that
$\soundness\,p$ reflects $p$.

Then, since $\completeness$ returns normal forms, we get that the composite function
$\completeness\,(\soundness\,p)$ evaluates to a normal form $q$
such that $q =_\beta p$.

Let us conclude this section by saying that the extension of this
proof to universal quantification and falsity, using so-called exploding nodes,
has been studied e.g. in~\cite{HerbelinLee09}. The extension to
first-order classical logic has been studied
e.g. in~\cite{IlikLeeHerbelin10}.
The case of disjunction and existential quantification is typically
addressed using variants of Kripke semantics~\cite{Ilik11b},
Beth models, topological models~\cite{CoquandSmith95}, or various alternative
semantics~(e.g. \cite{AbelSattler19,AltDybHofSco01,Okada02,Sambin95}).

One of the purposes of this paper is precisely to start a comparative
exploration of the computational contents of proofs of G\"odel's
completeness theorem and of the question of whether they provide a
normalisation procedure. In the case of Henkin's proof, the answer is
negative: even if the resulting object-level proof that will be
constructively obtained in Section~\ref{sec:Henkin} is related to the proof
of validity in the meta-logic, it is neither cut-free nor isomorphic
to it. In particular, it drops information from the validity proof by
sharing subparts that prove the same subformula as will be emphasised
in Section~\ref{sec:Henkin-S1}.

\subsection{The intuitionistic provability of the different statements of
completeness}
\label{sec:intuitionistic-provability}

Statements C1, C2 and C3, as well as statements S1, S2, S3 are
classically equivalent but not intuitionistically equivalent. In particular,
only C2 and S2 are intuitionistically provable.

More precisely, since our object language has only negative connectives,
the formula $\truth{\mathcal{T}}{\sigma}{\mathcal{M}}$ is in turn composed of
only negative connectives in the metalanguage. Hence, the only positive
connective in the statements C2 and S2 is the existential quantifier
asserting the existence of a model.

This existential quantifier is intuitionistically provable as our formulation of Henkin's
proof of S2 given in the next section shows: given a proof of
consistency of a theory, we can produce a syntactic model in the form
of a specific predicate. It shall however be noted that this predicate is not itself recursive
in general, since constructing this model is in general equivalent to
producing an infinite path in any arbitrary infinite binary tree
(such an infinite path is a priori not recursive, see
Kleene~\cite{Kleene52}, Simpson~\cite{Simpson99}).

Otherwise, from an intuitionistic point of view, statements C1 and S1 are particularly
interesting, as they promise to produce (object) proofs in
the object language out of proofs of validity in the
metalanguage. However, Kreisel~\cite{Kreisel58b} showed, using a result
by G\"odel~\cite{Godel31}, that C1 is equivalent to Markov's principle
over intuitionistic second-order arithmetic. This has been studied in
depth by McCarty~\cite{McCarty08} ant it turns out
that S1 is actually equivalent to Markov's principle if the theory is
recursively enumerable.

However, for arbitrary theories, reasoning
by contradiction on formulae of arbitrarily large logical complexity is
in general needed as the following adaptation of McCarty's proof
shows: Let $\mathcal{A}$ be an arbitrary formula of the metalanguage and consider e.g. the theory
defined by $B \in \mathcal{T} \defeq (B = \dbot) \land \mathcal{A} \lor (B = X) \land \neg \mathcal{A}$ for $X$ a distinct propositional atom of the object language.
We intuitionistically have that $\mathcal{T} \vDash X$ because
this is a negative formulation of a classically provable
statement\footnote{Let $\mathcal{M}$ and $\sigma$ such that
$\truth{\mathcal{T}}{\sigma}{\mathcal{M}}$. We first show $\neg \mathcal{A}$. Indeed, if $\mathcal{A}$
  holds, then $\dbot \in \mathcal{T}$ and we get by
$\truth{\mathcal{T}}{\sigma}{\mathcal{M}}$ that the model is contradictory. But if
  $\neg \mathcal{A}$, then $X \in \mathcal{T}$, hence $\truth{X}{\sigma}{\mathcal{M}}$.}.
By completeness, we get $\mathcal{T} \vdash X$,
and, by case analysis on the normal form of the so-obtained proof, we
infer that either $\mathcal{A}$ or $\neg \mathcal{A}$.

The need for Markov's principle is connected to how $\dbot$ is interpreted
in the model. Krivine~\cite{Krivine96b} showed that for a language without
$\dbot$\footnote{So-called minimal classical logic in~\cite{AriHer03}, which
  is however not functionally complete since no formula can represent the always-false function.},
C1 is provable intuitionistically. As
analysed by Berardi~\cite{Berardi99} and Berardi and
Valentini~\cite{BerardiValentini04}, Markov's principle is not needed
anymore if we additionally accept the degenerate model where all formulae
including $\dbot$ are interpreted as true\footnote{This is
  similar to the approach followed by Friedman~\cite{Friedman75} and
  Veldman~\cite{Veldman76} to intuitionistically prove the
  completeness of intuitionistic logic with respect to a relaxing of
  Beth models with so-called fallible models, and to a relaxing of
  Kripke models with so-called exploding nodes, respectively.}.
Let us formalise this precisely.

We define a {\em possibly-exploding} model $\mathcal{M}$ to be a
model\footnote{Such model is also called \emph{intuitionistic structure} in Krivtsov~\cite{Krivtsov10}.}
$(\mathcal{D},\mathcal{F},\mathcal{P},\mathcal{A}_{\dbot})$ such that $(\mathcal{D},\mathcal{F},\mathcal{P})$ is
a model in the previous sense and $\mathcal{A}_{\dbot}$ a fixed formula
intended to interpret $\dbot$ . The definition of truth is
then modified as follows:
$$
\begin{array}{lll}
\ptruth{\dbot}{\sigma}{\mathcal{M}} & \defeq & \mathcal{A}_{\dbot}\\
\end{array}
$$ with the rest of clauses unchanged\footnote{In~\cite{BerardiValentini04}, a
classical possibly-exploding model is called a minimal model, in reference to minimal logic~\cite{Johansson37}. The difference between non-exploding models and possibly-exploding models can actually be interpreted from the point of view of linear logic not as a difference of definition of models but as a difference of interpretations of the false connective. An non-exploding model is a model where the false connective is interpreted as the connective $0$ of linear logic (a positive connective, neutral for the standard disjunction and with no introduction rule). A possibly-exploding model is a model where the false connective is interpreted as the connective $\bot$ of linear logic (a negative connective informally standing for an empty sequent). See e.g. Okada~\cite{Okada02} or Sambin~\cite{Sambin95} for examples of differences of interpretation of $0$ and $\bot$ in completeness proofs for linear logic.}

Note that because $\dbot \dimp A$ is a
consequence of $\dneg\dneg A \dimp A$, the following holds
for all $A$ and all $\sigma$ in any
classical possibly-exploding model:
$$\mathcal{A}_{\dbot} ~\imp~ \ptruth{A}{\sigma}{\mathcal{M}}\,,$$ so we do not need to further enforce
it\footnote{As a matter of purity, since it is standard that the
  classical scheme $\dneg\dneg A \dimp A$ is equivalent to the
  conjunction of a purely classical part, namely Peirce's law
  representing the scheme $((A\dimp B)\dimp A)\dimp B$ and of a purely
  intuitionistic part, namely ex falso quodlibet representing the
  scheme $\dbot \dimp A$, we could have decomposed ${\classic}$ into
  the union of $\peirce \defeq \{((A \dimp B) \dimp A) \dimp A
  ~|~ A \in\Form\}$ and of $\exfalso \defeq \{\dbot\dimp A ~|~ A
  \in\Form\}$.

  As already said in Section~\ref{sec:models}, the condition
  $\truth{\classic}{\sigma}{\mathcal{M}}$, and in particular the conditions
  $\truth{\peirce}{\sigma}{\mathcal{M}}$ and $\truth{\exfalso}{\sigma}{\mathcal{M}}$
  are needed to show soundness with respect to classical models
  in a minimal setting.  In an intuitionistic setting,
  $\truth{\exfalso}{\sigma}{\mathcal{M}}$ holds by default and does not
  have to be explicitly enforced. In a classical setting,
  $\truth{\peirce}{\sigma}{\mathcal{M}}$ does not have to be explicitly
  enforced. So, requiring these conditions is to ensure that the
  definition of validity is the one we want independently of the
  specific properties of the metalanguage.

  In contrast, for the purpose of completeness, possibly-exploding models
  are needed for an intuitionistic proof of C1 to be possible, but none
  of $\peirce(\mathcal{M})$ and $\exfalso(\mathcal{M})$ are required.}.  Let
us rephrase C1 and S1 using classical possibly-exploding models:

$$
\begin{array}{ll}
C1'. & \pvDash A ~\imp~\isprovable{}{A}\\
S1'. & \mathcal{T} \pvDash A ~\imp~\isprovable{\mathcal{T}}{A}\\
\end{array}
$$
where
$$
\mathcal{T} \pvDash A ~ \defeq ~ \forall \mathcal{M}\, \forall \sigma\, (\ptruth{\classic}{\sigma}{\mathcal{M}} \imp \ptruth{\mathcal{T}}{\sigma}{\mathcal{M}} \imp \ptruth{A}{\sigma}{\mathcal{M}})\\
\\
$$

In particular, it is worthwhile to notice that $\mathcal{T} \pvDash A$
and $\mathcal{T} \vDash A$ are classically equivalent since, up to logical equivalence, $\pvDash$
only differs from $\vDash$ by an extra quantification over the
degenerate always-true model. Hence C1 and C1',
as well as S1 and S1', are classically equivalent too. But C1' as well
as S1' are intuitionistically provable for recursively enumerable
theories, while C1 and S1, even for recursively enumerable
theories, would require Markov's principle\footnote{Markov's principle
  can actually be ``intuitionistically'' implemented e.g. by using an
  exception mechanism~\cite{Herbelin10}, so a computational content to
  weak completeness and strong completeness for recursively enumerable
  theories could as well be obtained without any change to the interpretation of
  $\dbot$.}.

Let us conclude this section by considering C3 and
S3. These statements are not intuitionistically provable:
if they were, provability could be decided. This does not mean however
that we cannot compute with C3 and S3. Classical logic is
computational (see e.g.~\cite{Parigot92}), but for an evaluation to
be possible, an interaction with a proof of a statement which uses
C3 or S3 is needed. The formalisation by Blanchette, Popescu and
Traytel~\cite{BlanchettePopescuTraytel14} might be a good starting
point to analyse the proof but we will not explore this further here.

\subsection{Chronology and recent related works}

The extraction of a computational content from Henkin's proof was
obtained by the authors from an analysis of the
formalisation~\cite{Ilik08} in the Coq proof assistant~\cite{Coq21}
of Henkin's proof. It was first presented at the TYPES
conference in Warsaw, 2010. The paper was essentially written in 2013
but it remained in draft and unstable form until 2016. A
non-peer-reviewed version was made publicly available on the web page
of the first author late 2016 and slightly updated in 2019. Our
constructive presentation of Henkin's proof inspired Forster, Kirst
and Wehr to formalise the completeness theorem in the intuitionistic
setting of Coq~\cite{ForsterKirstWehr20,ForsterKirstWehr21}. This
encouraged us in 2022 to polish the paper one step further with the
objective of a submission.

A new version was made public on the web page of the first author in
November 2022 together with an
implementation\footnote{\url{http://herbelin.xyz/articles/henkin.ml}}
of the proof in the OCaml programming language. The new version
presented the proof in a more structured way and
made explicit the need for a weakening rule left unadressed in the
2016 version. It also provided an extension of the proof to the
disjunction conjunctive, relying on (a slight generalisation of) the
weakly classical principle of Double Negation Shift ($\DNS$). Finally,
the current and last version includes the extra result that completeness in
the presence of disjunction requires at least a variant of a weak form
of Double Negation Shift that Kirst recently
identified~\cite{Kirst23}.

Various investigations of G\"odel's completeness in an intuitionistic
setting have been published since our proof was first presented. To
our knowledge, in addition
to~\cite{ForsterKirstWehr20,ForsterKirstWehr21,Kirst23},
this includes papers by Krivtsov~\cite{Krivtsov15} in intuitionistic
arithmetic and
Espíndola~\cite{EspindolaPhD,Espindola16} in intuitionistic set
theory (IZF). Their results are discussed in the
section~\ref{sec:disjunction}.

\section{The computational content of Henkin's proof of Gödel's completeness}
\label{sec:Henkin}

We now recall Henkin's proof of completeness and analyse its
computational content.

\subsection{Henkin's proof, slightly simplified}

We give a simplified form of Henkin's proof of the strong
form of Gödel's completeness theorem~\cite{Henkin49a}, formulated as a model existence theorem, that is as
statement S2. The simplification is on the use of free variables
instead of constants in Henkin axioms and in the use of only
implicative formulae in the process of completion of a consistent set
of formulae into a maximally consistent
one\footnote{Smullyan~\cite{Smullyan68} credits
  Hasenjaeger~\cite{Hasenjaeger53} and Henkin independently for proof variants using
  free variables (see also Henkin~\cite{Henkin96}). In particular,
  this allows to build a maximal consistent theory in one step instead
  of a countable number of steps as in Henkin's original proof. See
  also~\cite[Th. IV.3.3]{Simpson99} for a proof building a maximal
  consistent theory in one step.}.

Let $\mathcal{T}$ be a consistent set of formulae mentioning an at most
countable\footnote{In the presence of uncountably many symbols, one
  would need the ultrafilter theorem to build the model and this would
  require extra computational tools to make the proof
  constructive. See~\cite{Jech73} for the equivalence in set theory
  between the ultrafilter theorem and the completeness theorem on
  non-necessarily countable signatures.}
number of function symbols and predicate symbols. Let $\mathcal{X}_1$ and
$\mathcal{X}_2$ be countable sets of variables forming a partition of
$\mathcal{X}$. We can assume without loss of generality that the free
variables of the formulae in $\mathcal{T}$ are in $\mathcal{X}_1$ leaving
$\mathcal{X}_2$ as a pool of fresh variables.

We want to show that $\mathcal{T}$ has a model, and for that purpose, we
shall complete it into a consistent set $\mathcal{S}_{\omega}$ of
formulae which is maximal in the sense that if $A \not\in \mathcal{S}_{\omega}$
then $\dneg A \in \mathcal{S}_{\omega}$. We shall also
ensure that for every universally quantified formula $\dforall
x\,A(x)$, there is a corresponding so-called Henkin axiom $A(y) \dimp
\dforall x\,A(x)$ in $\mathcal{S}_{\omega}$ with $y$ fresh in $\dforall
x\, A(x)$. For the purpose of this construction, we fix an injective
enumeration $\phi$ of formulae of the form $\dforall x\,A(x)$ or
$A\dimp B$ and write $\godel{A}$ for the index of a formula $A$ of such form in
the enumeration. We also take $\phi$ so that formulae of even index
are of the form $\dforall x\,A(x)$ and formulae of odd index are of
the form $A\dimp B$.

Let $\mathcal{S}_0$ be $\mathcal{T}$ and assume that we have already built
$\mathcal{S}_n$. If $n$ is even, $\phi(n)$ has the form $\dforall
x\,A(x)$. We then consider a variable $x_{n/2}\in \mathcal{X}_2$ which is fresh in
all $\phi(i)$ for $i\leq n$ and we set $\mathcal{S}_{n+1} \defeq
\mathcal{S}_n \cup (A(x_{n/2}) \dimp \dforall x\,A(x))$. Otherwise, $\phi(n)$ is
an implicative formula and we consider two cases. If $\mathcal{S}_n \cup
\phi(n)$ is consistent, i.e., if $(\mathcal{S}_n \cup \phi(n) \vdash
\dbot) \imp \bot$, we set $\mathcal{S}_{n+1} \defeq \mathcal{S}_n \cup
\phi(n)$. Otherwise, we set $\mathcal{S}_{n+1} \defeq \mathcal{S}_n$.
We finally define the predicate $A\in \mathcal{S}_{\omega} \defeq \exists
n\,(\mathcal{S}_n \vdash A)$, i.e. $\exists n\exists \Gamma \subset
\mathcal{S}_n\; (\Gamma \vdash A)$, and this is the basis of a syntactic model
$\mathcal{M}_0$ defined by taking

$$
\begin{array}{lll}
\mathcal{D} & \defeq & \Term\\
\mathcal{F}(f)(t_1,...,t_n) & \defeq & f(t_1,...,t_n)\\
\mathcal{P}(P)(t_1,...,t_n) & \defeq & P(t_1,...,t_n) \in \mathcal{S}_{\omega}\\
\end{array}
$$

By induction, each $\mathcal{S}_n$ is consistent. Indeed, if $\phi(n)$ is
implicative and $\mathcal{S}_{n+1} \equiv \mathcal{S}_n \cup \phi(n)$, it is
precisely because $\mathcal{S}_{n+1}$ is consistent. Otherwise, the
consistency of $\mathcal{S}_{n+1}$ comes from the consistency of $\mathcal{S}_n$.
If $\phi(n)$ is some $\dforall x\,A(x)$, then $\mathcal{S}_{n+1}
\equiv \mathcal{S}_n \cup (A(x_{n/2}) \dimp \dforall x\,A(x))$. This is
consistent by freshness of $x_{n/2}$ in both $\mathcal{T}$ and in the
$\phi(i)$ for $i\leq n$. Indeed, because $x_{n/2}$ is fresh, any
proof of $\mathcal{S}_n \cup (A(x_{n/2}) \dimp \dforall x\,A(x)) \vdash \dbot
$ can be turned into a proof of $\mathcal{S}_n \cup \dneg \dforall y\, \dneg (A(y)
\dimp \dforall x\,A(x)) \vdash \dbot $, which itself can be turned
into a proof of $\mathcal{S}_n \vdash \dbot $ since $\dneg \dforall y\, \dneg \,(A(y)
\dimp \dforall x\,A(x))$ is a classical tautology\footnote{The famous ``Drinker Paradox''.}

Let $A$ be a formula and $\sigma$ a substitution of its free variables. We now show by induction on
the logical depth\footnote{In particular, we consider $B(t)$ to be smaller to
  $\forall x\,B(x)$ for any $t$.} of $A$ that $\truth{A}{\sigma}{\mathcal{M}_0}
  \myiff A[\sigma] \in \mathcal{S}_{\omega}$, where $A[\sigma]$ denotes the result of substituting the free variables of $A$ by the terms in $\sigma$. This is sometimes considered an easy
combinatoric argument but we shall detail the proof because it is here
that the computational content of the proof is non-trivial. Moreover,
we do not closely follow Henkin's proof who is making strong use of
classical reasoning. We shall instead reason intuitionistically, which
does not raise any practical difficulty here.

\begin{itemize}

\item Let us focus first on the case when $A$ is $B\dimp C$. One way to
  show $B[\sigma] \dimp C[\sigma] \in \mathcal{S}_{\omega}$ from $\truth{B\dimp
    C}{\sigma}{\mathcal{M}_0}$ is to show that for $n$ being $\godel{B[\sigma]\dimp
    C[\sigma]}$, the set $\mathcal{S}_n \cup B[\sigma] \dimp C[\sigma]$ is consistent,
  i.e. that a contradiction arises from $\mathcal{S}_n \cup B[\sigma] \dimp C[\sigma]
  \vdash \dbot$. Indeed, from the latter, we get both $\mathcal{S}_n
  \vdash B[\sigma]$ and $\mathcal{S}_n \vdash \dneg C[\sigma]$. From $\mathcal{S}_n \vdash
  B[\sigma]$ we get $\truth{B}{\sigma}{\mathcal{M}_0}$ by induction hypothesis, hence
  $\truth{C}{\sigma}{\mathcal{M}_0}$ by assumption on the truth of $B\dimp
  C$. Then $C[\sigma] \in \mathcal{S}_{\omega}$ again by induction hypothesis,
  hence $\mathcal{S}_{n'} \vdash C[\sigma]$ for some $n'$. But also $\mathcal{S}_n
  \vdash \dneg C[\sigma]$, hence $\mathcal{S}_{\mathit{max}(n,n')} \vdash \dbot$ which
  contradicts the consistency of $\mathcal{S}_{\mathit{max}(n,n')}$.

\item Conversely, if $B[\sigma] \dimp C[\sigma] \in \mathcal{S}_{\omega}$, this means
  $\mathcal{S}_n \vdash B[\sigma] \dimp C[\sigma]$ for some $n$. To prove $\truth{B\dimp
  C}{\sigma}{\mathcal{M}_0}$, let us assume $\truth{B}{\sigma}{\mathcal{M}_0}$. By
  induction hypothesis we get $\mathcal{S}_{n'} \vdash B[\sigma]$ for some $n'$ and
  hence $\mathcal{S}_{\mathit{max}(n,n')} \vdash C[\sigma]$, i.e. $C[\sigma] \in \mathcal{S}_{\omega}$.
  We conclude by induction hypothesis to get
  $\truth{C}{\sigma}{\mathcal{M}_0}$.

\item Let us then focus on the case when $A$ is $\dforall
  x\,B$. For $n$ even being $\godel{(\dforall x\,B)[\sigma]}$, we have
  $(B[\sigma,x\leftarrow x_{n/2}] \dimp (\dforall x\,B)[\sigma]) \in \mathcal{S}_{n+1}$. From
  $\truth{\dforall x\,B(x)}{\sigma}{\mathcal{M}_0}$ we get
  $\truth{B(x)}{\sigma, x\leftarrow x_{n/2}}{\mathcal{M}_0}$ and by the induction hypothesis
    we then get the existence of some $n'$ such that $\mathcal{S}_{n'}
    \vdash B[\sigma,x \leftarrow x_{n/2}]$. Hence, $\mathcal{S}_{\mathit{max}(n+1,n')} \vdash (\dforall
    x\,B)[\sigma]$, which means $(\dforall x\,B)[\sigma] \in \mathcal{S}_{\omega}$.

\item Conversely, assume $\mathcal{S}_n \vdash (\dforall x\,B)[\sigma]$ for some
  $n$ and prove $\truth{\forall x\,B}{\sigma}{\mathcal{M}_0}$. Let $t$ be
  a term.
  From $\mathcal{S}_n \vdash (\dforall x\,B)[\sigma]$ we get $\mathcal{S}_n \vdash
  B[\sigma,x\leftarrow t]$, hence, by induction hypothesis,
  $\truth{B(x)}{\sigma,x\leftarrow t}{\mathcal{M}_0}$.

\item Let us then consider the case $A$ is $\dbot$. By ex falso
  quodlibet in the metalanguage, it is direct that $\bot \imp (\dbot \in \mathcal{S}_{\omega})$.

\item Conversely, let us prove $(\dbot \in \mathcal{S}_{\omega}) \imp
  \bot$. From $\dbot \in \mathcal{S}_{\omega}$ we know $\mathcal{S}_n \vdash
  \dbot$ for some $n$ which, again, contradicts the consistency of
  $\mathcal{S}_n$.

\item The case when $A$ is $P(t_1,\ldots,t_n)$ is by definition.
\end{itemize}

Before completing the proof, it remains to prove that the model is
classical. Using the equivalence between $\truth{A}{\substid}{\mathcal{M}_0}$
and $A \in \mathcal{S}_{\omega}$ for $A$ closed and $\substid$
the empty substitution, it is enough to prove that $\dneg\dneg
A\in \mathcal{S}_{\omega}$ implies $A\in \mathcal{S}_{\omega}$. But the
former means $\mathcal{S}_{n} \vdash \dneg\dneg A$ for some $n$, hence
$\mathcal{S}_{n} \vdash A$ by classical reasoning in the object language,
hence $A \in \mathcal{S}_{\omega}$.

We are now ready to complete the proof: for every $B\in \mathcal{T}$,
since $\mathcal{T} \vdash B$, we get $B \in \mathcal{S}_\omega$ and hence
$\truth{B}{\substid}{\mathcal{M}_0}$.

\subsection{From Henkin's proof to a proof with respect to possibly-exploding models}
\label{sec:Henkin-S1}

Let us fix a formula $A_0$ and a recursively enumerable theory
$\mathcal{T}_0$, i.e. a theory defined by a $\Sigma^0_1$-statement.  To
get a proof of statement S1 for $\mathcal{T}_0$ and $A_0$ is easy by
using Markov's principle: to prove $\mathcal{T}_0 \vdash A_0$ from $\mathcal{T}_0
\vDash A_0$, let us assume the contrary, namely that
$\mathcal{T}_0 \cup \dneg A_0$ is consistent. Then, we can complete
$\mathcal{S}_0 \defeq \mathcal{T}_0 \cup \dneg A_0$ into $\mathcal{S}_{\omega}$ and
build out of it a classical model $\mathcal{M}_0$ such that $\forall
B\in\mathcal{T}_0\, \truth{B}{\substid}{\mathcal{M}_0}$ as well as $\truth{\dneg
  A_0}{\substid}{\mathcal{M}_0}$, i.e.  $\neg (\truth{A_0}{\substid}{\mathcal{M}_0})$. But
this contradicts $\mathcal{T}_0 \vDash A_0$ and, because $\mathcal{T}_0$ is
$\Sigma^0_1$, hence $\mathcal{T}_0 \vdash A$ as well, Markov's principle applies.

As discussed in Section~\ref{sec:intuitionistic-provability}, S1
cannot be proved without Markov's principle, so we shall instead prove S1'. To
turn the proof of S2 into a proof of S1' which does not require
reasoning by contradiction, we shall slightly change the construction of
$\mathcal{S}_{\omega}$ from $\mathcal{T}_0 \cup \dneg A_0$ so that it is not consistent in an
absolute sense, but instead consistent relative to $\mathcal{T}_0 \cup \dneg A_0 $. In
particular, we change the condition for extending $\mathcal{S}_{2n+1}$
with $\phi(2n+1)$ to be that $\mathcal{S}_{2n+1} \cup \phi(2n+1)$ is {\em
  consistent relative} to $\mathcal{T}_0 \cup \dneg A_0$.

Then, we show by induction not that $\mathcal{S}_n$ is consistent but
that its inconsistency reduces to the inconsistency of $\mathcal{T}_0 \cup \dneg A_0$.

For the construction of the now possibly-exploding model, we take as
interpretation of $\dbot$ the formula $\mathcal{T}_0, \dneg A_0 \vdash \dbot$.
Proving $\dbot \in \mathcal{S}_{\omega} \imp \ptruth{\dbot}{\sigma}{\mathcal{M}_0}$
now reduces to proving $\mathcal{S}_{n} \vdash \dbot \imp
\mathcal{T}_0, \dneg A_0 \vdash \dbot$ which is the statement of relative
consistency\footnote{Interestingly enough, since $\mathcal{T}_0 \vdash
  A$ effectively holds as soon as an effective proof of validity of $A$ is given, the model
  we build is then the degenerate one in which all formulae are
  true.}. Conversely, $\ptruth{\dbot}{\sigma}{\mathcal{M}_0} \imp \dbot \in
\mathcal{S}_{\omega}$ now comes by definition of $\mathcal{S}_{0}$.

The change in the definition of $\mathcal{S}_\omega$ as well as the use
of possibly-exploding models is connected to Friedman's
$A$-translation~\cite{Friedman78} being able to turn Markov's
principle into an admissible rule. Here, $A$ is the $\Sigma^0_1$-formula $\mathcal{T}_0 \vdash
A_0$ and by replacing $\bot$ with $A$ in the definition of model, hence of validity, as
well as in the definition of $\mathcal{S}_\omega$, we are able to prove
$(A\imp A)\imp A$ whereas only $(A\imp \bot)\imp \bot$ was otherwise
provable. Then, $A$ comes trivially from $(A\imp A)\imp
A$.

This was the idea followed by Krivine~\cite{Krivine96b} in his
constructive proof of G\"odel's theorem for a language restricted to
$\dimp$ and $\dforall$, as analysed and clarified in Berardi and
Valentini~\cite{BerardiValentini04}.

As a final remark, one could wonder whether the construction of $\mathcal{S}_{2n+2}$
by case on an undecidable statement is compatible with
intuitionistic reasoning. Indeed, constructing the sequence of
formulae added to $\mathcal{T}_0 \cup \dneg A_0$ in order to get $\mathcal{S}_n$ seems to
require a use of excluded-middle. However, in the proof of
completeness, only the property $A\in \mathcal{S}_n$ matters, and this
property is directly definable by induction as

$$
\begin{array}{lll}
A\in \mathcal{S}_{0} & \defeq & A \in \mathcal{T}_0 \cup \dneg A_0\\
A \in \mathcal{S}_{n+1} & \defeq & A \in \mathcal{S}_n\\
& & \lor~ (\exists p\,(n = 2p+1) ~\land~ \phi(n) = \dforall x\,B(x) ~\land~ A =
(B(x_{p+1}) \,\dimp\, \dforall x\,B(x))\\
& & \lor~ (\exists p\,(n = 2p) ~\land~
(\mathcal{S}_n, A \vdash \dbot \imp \mathcal{T}_0, \dneg A_0 \vdash \dbot)  ~\land~ A = \phi(n))
\end{array}
$$

Note however that $\mathcal{S}_n$ is used in negative position of an
implication in the definition of $\mathcal{S}_{n+1}$. Hence, the
complexity of the formula $A \in \mathcal{S}_n$ seen as a type of
functions is a type of higher-order functions of depth
$n$.

\subsection{The computational content of the proof of completeness}

We are now ready to formulate the proof as a program. We shall place
ourselves in an axiom-free second-order intuitionistic arithmetic
equipped with a proof-as-program interpretation\footnote{A typical
  effective framework for that purpose would be a fragment of the
  Calculus of Inductive Constructions such as it is implemented in the
  Coq proof assistant~\cite{Coq19a}, or Matita~\cite{Matita}. The
  Calculus of Inductive Constructions is an impredicative extension of
  Martin-L\"of's type theory~\cite{MartinLof75}.}, as already
considered in Section~\ref{sec:semantic-normalisation}. Additionally,
we shall identify the construction of existentially quantified
formulae and the construction of proofs of conjunctive formulae. For
instance, we shall use the notation $(p_1,...,p_n)$ for the proof of an $n$-ary
combination of existential quantifiers and conjunctions. We shall also write
$\dest{p}{(x_1,...,x_n)}{q}$ for a proof obtained by decomposition of
the proof $p$ of an $n$-ary combination of existential quantifiers and conjunctions. We
shall write $\efq\,p$ for a proof of $A$ obtained by ex falso quodlibet from a proof $p$ of $\bot$.

We shall use the letters $n$, $A$, $\Gamma$,
$m$, $p$, $q$, $r$, $h$, $g$, $f$, $k$ and their variants to refer to natural
numbers, formulae, contexts of formulae, proofs of truth, proofs of
derivability in the object language, proofs of belonging to $\mathcal{S}_\omega$,
proofs of inconsistency from adding an implicative formula to $\mathcal{S}_{2n+1}$,
proofs of belonging to $\mathcal{T}_0$,
proofs of inclusion in $\mathcal{T}_0$, proofs of inclusion in
extensions of $\mathcal{T}_0$, proofs of relative consistency,
respectively.

The key property is $A\in \mathcal{S}_{\omega}$ which unfolds as $\exists
n\,\exists \Gamma\, (\Gamma \subset \mathcal{S}_n\wedge \Gamma \vdash
A)$. Rather than defining $\Gamma \subset \mathcal{S}_n$ from $A \in
\mathcal{S}_n$ and the latter by cases, we now directly take $\Gamma \subset
\mathcal{S}_n$ as our primitive concept, so that defining $A \in \mathcal{S}_n$
is actually not needed anymore. Rephrasing the property that
$\mathcal{S}_{n}$ is inconsistent in terms of $\Gamma \subset
\mathcal{S}_n$ is easy: it is enough to tell that $\Gamma \vdash \dbot$
for some $\Gamma \subset \mathcal{S}_{n}$. In particular, we can
write $\mathcal{S}_n, A \vdash \dbot$ to mean $\exists
\Gamma\,(\Gamma \subset \mathcal{S}_{n}\;\land\; \Gamma, A \vdash
\dbot)$.

We first define by cases the predicate $\Gamma \subset
\mathcal{T}_0$:

$$
\seqr{\Jbase}{}{\epsilon \,\subset\, \mathcal{T}_0}
\qquad
\qquad
\seqr{\Jcons}{\Gamma \subset \mathcal{T}_0\qquad A\in\mathcal{T}_0}{\Gamma, A\, \subset\, \mathcal{T}_0}
$$

Then, we can define $\Gamma \subset \mathcal{S}_n$ by cases: a formula
$B$ is allowed to occur in such $\Gamma$ either because it is in
$\mathcal{T}_0$ (clause $J_{\mathit{cons}}$), or because it is an Henkin
axiom added at step $2n$ (clause $I_{\forall}$), or because it is an
implication added at step $2n+1$ together with a proof of relative
consistency of $\mathcal{S}_{n+1} \cup B$ with respect to $\mathcal{T}_0
\cup \dneg A_0$ (clause $I_{\imp}$). Note that we enforce in all cases
that at least $\dneg A_0$ is in such $\Gamma$ (clause $\IO$) and that
we can always skip adding a formula at some step of the construction
(clause $\IS$). Formally, the definition is:

$$
\seqr{\IO}{\Gamma \subset \mathcal{T}_0}{\Gamma, \dneg A_0 \subset \mathcal{S}_0}
\qquad
\qquad
\seqr{\IS}{\Gamma \subset \mathcal{S}_n}{\Gamma \subset \mathcal{S}_{n+1}}
\qquad
\qquad
\seqr{\Iforall}{\Gamma \subset \mathcal{S}_{2n}}{\Gamma,\; A(x_n) \,\dimp\, \dforall x\,A(x) \,\subset\, \mathcal{S}_{2n+1}}
$$

$$
\seqr{\Iimp}{\Gamma \subset \mathcal{S}_{2n+1} \qquad \mathcal{S}_{2n+1}, A\dimp B \vdash \dbot \imp \mathcal{T}_0, \dneg A_0 \vdash \dbot}{\Gamma,\, A\,\dimp\, B \,\subset\, \mathcal{S}_{2n+2}}
$$
where $\phi(2n)$ is $\dforall x\,A(x)$ in $I_{\forall}$ and
$\phi(2n+1)$ is $A\dimp B$ in $I_{\imp}$.

We now write as a program the proof that $\mathcal{S}_\omega$ is
consistent relative to $\mathcal{T}_0 \cup \dneg A_0$, i.e. that
$\dbot \in \mathcal{S}_{\omega}$ implies $\mathcal{T}_0, \dneg A_0 \vdash \dbot$. The latter
expands to $\exists \Gamma\,(\Gamma\subset \mathcal{T}_0\;\land\; \Gamma,
\dneg A_0 \vdash \dbot)$ which we see as made of triples of
the form $(\Gamma,g,p)$, with $\Gamma$ a context, $g$ a proof of
$\Gamma\subset \mathcal{T}_0$ and $p$ a proof of $\Gamma', \dneg A_0
\vdash \dbot$.

This proof, which we call $\squeeze_n^\Gamma$, takes as arguments
a quadruple $(n,\Gamma,f,p)$ where $f$ is a proof of $\Gamma \subset
  \mathcal{S}_n$ and $p$ proof of $\Gamma \vdash \dbot$. It proceeds by cases on the
  proof of $\Gamma \subset \mathcal{S}_n$. When extended at odd $n$, it works by
calling the \emph{continuation} justifying that adding the formula
$\phi(2p+1)$ preserves consistency, and, when extended at even $n$, by composing the
resulting proof of inconsistency with a proof of the Drinker's paradox
($\drinker_y$ is the proof which builds a proof of $\Gamma \vdash
\dbot$ from a proof of $\Gamma, A(y) \dimp \dforall x\,A(x) \vdash
\dbot$, knowing that $y$ does not occur in $\Gamma, \dforall x\, A(x)$,
see Figure~\ref{fig:rules}).
$$
\begin{array}{lllll}
\squeeze & {:} & \dbot \in \mathcal{S}_{\omega} & {\imp} & \mathcal{T}_0, \dneg A_0 \vdash \dbot\\
\squeeze & & (0,(\Gamma,\dneg A_0),\IO\,g,p) & \defeq & (\Gamma,g,p)\\
\squeeze & & (n+1,\Gamma,\IS\,f,p) & \defeq & \squeeze\, (n,\Gamma,f,p)\\
\squeeze & & (2n+1,(\Gamma,A),\mathtt{I_\forall}\,f,p) & \defeq & \squeeze\, (2n,\Gamma,f,\drinker_{x_n}\,p)\\
\squeeze & & (2n+2,(\Gamma,A),\Iimp\,(f,k),p) & \defeq & k\;(\Gamma,f,p)\\
\end{array}
$$

A trivial lemma implicit in the natural language formulation of the
proof of completeness is the lemma asserting $\dneg A_0 \subset
\mathcal{S}_n$. The proof is by induction on $n$:
$$
\begin{array}{lclllcll}
\inj_n & {:} & {\dneg A_0 \subset \mathcal{S}_n} \quad~~ \\
\inj_0 & \defeq & \IO(J_{\mathit{base}}) \\
\inj_{n+1} & \defeq & \IS (\inj_n)\\
\end{array}
$$

A boring lemma which is implicit in the proof of completeness in natural
language is that $\Gamma \subset \mathcal{S}_n$ and $\Gamma'\subset \mathcal{S}_{n'}$
imply $\Gamma \cup \Gamma' \subset \mathcal{S}_{\mathit{max}(n,n')}$.
It looks obvious because one tends to think of
$\Gamma \subset \mathcal{S}_n$ as denoting the inclusion of $\Gamma$
within a uniquely defined relatively consistent set $\mathcal{S}_n$.
However, the computational approach to the proof shows that
$\mathcal{S}_n$ has no computational content per se: only proofs of $\Gamma \subset \mathcal{S}_n$ have, and such proofs are collections
of proofs of relative
consistency for only those implicative formulae which are in $\Gamma$.
These formulae are those
inspected by the lemma $A\in \mathcal{S}_{\omega} \Leftrightarrow
\ptruth{A}{\sigma}{\mathcal{M}_0}$, which in practice are subformulae of
the formulae in $\mathcal{T}_0$.

For $\Gamma_1$ and $\Gamma_2$ included in $\mathcal{T}_0$, we write
$\Gamma_1, \Gamma_2$ for their concatenation (possibly with redundancies). Otherwise, by
construction, any $\Gamma$ included in $\mathcal{S}_n$ for some $n$ has
either the form $\Gamma', \dneg A_0$ where $\Gamma'$ is included in
$\mathcal{T}_0$, or the form $\Gamma', A$ where $A$ has been added in
the process of enumeration.  We can then define $\Gamma_1 \cup
\Gamma_2$ for $\Gamma_1$ and $\Gamma_2$ included in
$\mathcal{S}_n$ for some $n$ by cases\footnote{Strictly speaking, this
  decomposition of any $\Gamma$ included in $\mathcal{S}_n$ for some $n$
  should be part of the structure of $\Gamma$ so as to be able to
  compute with it.}:

$$
\begin{array}{llllll}
\Gamma_1, \dneg A_0 &\cup& \Gamma_2, \dneg A_0 & ::= & \Gamma_1, \Gamma_2, \dneg A_0\\
\Gamma_1, \dneg A_0 &\cup& \Gamma_2, A & ::= & (\Gamma_1, \dneg A_0 ~\cup~ \Gamma_2),A\\
\Gamma_1,A &\cup& \Gamma_2, \dneg A_0 & ::= & (\Gamma_1 ~\cup~ \Gamma_2, \dneg A_0),A\\
\Gamma_1,A &\cup& \Gamma_2,A & ::= & (\Gamma_1 ~\cup~ \Gamma_2),A \quad & \\
\Gamma_1,A &\cup& \Gamma_2,B & ::= & (\Gamma_1,A ~\cup ~\Gamma_2),B \quad & \mbox{if $\godel{B} > \godel{A}$} \\
\Gamma_1,A &\cup& \Gamma_2,B & ::= & (\Gamma_1 ~\cup~ \Gamma_2,B),A & \mbox{if $\godel{A} > \godel{B}$} \\
\end{array}
$$

We can then define the merge of two proofs of
$\Gamma\subset\mathcal{S}_n$ by distinguishing when the contexts are
considered as subsets of $\mathcal{T}_0$ or as subsets of some
$\mathcal{S}_n$:

$$
\begin{array}{lrlcllcll}
{\subsetjoin^{\Gamma_1\Gamma_2}} : & \!\!\!\! & {\Gamma_1 \subset \mathcal{T}_0} &\!\!\!{\land}\!\!\!& {\Gamma_2 \subset \mathcal{T}_0} &&{\!\!\imp\!\!}& {\Gamma_1, \Gamma_2 \subset \mathcal{T}_0}\\
\subsetjoin^{\Gamma_1\epsilon} & (\!\!\!\! &g_1&\!\!\!,\!\!\!&\Jbase &\!\!\!\!)&\!\!\defeq\!\!& g_1 &\\
\subsetjoin^{\Gamma_1(\Gamma_2,A_2)} & (\!\!\!\! &g_1&\!\!\!,\!\!\!&\Jcons(g_2,h_2) &\!\!\!\!)&\!\!\defeq\!\!& \Jcons(\subsetjoin^{\Gamma_1\Gamma_2}(g_1,g_2),h_2) & \\
\\
{\join^{\Gamma_1\Gamma_2}_{n}} : & \!\!\!\! & {\Gamma_1 \subset \mathcal{S}_{n}} &\!\!\!{\land}\!\!\!& {\Gamma_2 \subset \mathcal{S}_{n}} &&{\!\!\imp\!\!}& {\Gamma_1 \cup \Gamma_2 \subset \mathcal{S}_n}\\
\join^{(\Gamma_1,\dneg A_0)(\Gamma_2,\dneg A_0)}_{0}\!\!\!\! & (\!\!\!\! &\IO(g_1)&\!\!\!,\!\!\!&\IO(g_2) &\!\!\!\!)&\!\!\defeq\!\!& \IO(\subsetjoin^{\Gamma_1\Gamma_2}(g_1,g_2) & \\
\join^{(\Gamma_1A)(\Gamma_2A)}_{2n+1} & (\!\!\!\! &\mathtt{I_\forall}(f_1)&\!\!\!,\!\!\!&\mathtt{I_\forall}(f_2) &\!\!\!\!)&\!\!\defeq\!\!& \mathtt{I_\forall}(\join^{\Gamma_1\Gamma_2}_{2n} (f_1, f_2)) & \\
\join^{(\Gamma_1A)(\Gamma_2A)}_{2n+2} & (\!\!\!\! &\mathtt{I_\imp}(f_1,k_1)&\!\!\!,\!\!\!&\mathtt{I_\imp}(f_2,k_2) &\!\!\!\!)&\!\!\defeq\!\!& \mathtt{I_\imp}(\join^{\Gamma_1\Gamma_2}_{2n+1} (f_1, f_2),k_1) & \\
\join^{(\Gamma_1A)\Gamma_2}_{2n+1} & (\!\!\!\! &\mathtt{I_\forall}(f_1)&\!\!\!,\!\!\!&\IS(f_2) &\!\!\!\!)&\!\!\defeq\!\!& \mathtt{I_\forall}(\join^{\Gamma_1\Gamma_2}_{2n} (f_1, f_2)) \\
\join^{(\Gamma_1A)\Gamma_2}_{2n+2} & (\!\!\!\! &\mathtt{I_\imp}(f_1,k_1)&\!\!\!,\!\!\!&\IS(f_2) &\!\!\!\!)&\!\!\defeq\!\!& \mathtt{I_\imp}(\join^{\Gamma_1\Gamma_2}_{2n+1} (f_1, f_2),k_1) \\
\join^{\Gamma_1(\Gamma_2A)}_{2n+1} & (\!\!\!\! &\IS(f_1)&\!\!\!,\!\!\!&\mathtt{I_\forall}(f_2) &\!\!\!\!)&\!\!\defeq\!\!& \mathtt{I_\forall}(\join^{\Gamma_1\Gamma_2}_{2n} (f_1, f_2)) &\\
\join^{\Gamma_1(\Gamma_2A)}_{2n+2} & (\!\!\!\! &\IS(f_1)&\!\!\!,\!\!\!&\mathtt{I_\imp}(f_2,k_2) &\!\!\!\!)&\!\!\defeq\!\!& \mathtt{I_\imp}(\join^{\Gamma_1\Gamma_2}_{2n+1} (f_1, f_2),k_2) &\\
\join^{\Gamma_1\Gamma_2}_{n+1} & (\!\!\!\! &\IS(f_1)&\!\!\!,\!\!\!&\IS(f_2) &\!\!\!\!)&\!\!\defeq\!\!& \IS(\join^{\Gamma_1\Gamma_2}_{n} (f_1, f_2)) &\\
\\
{\hjoin^{\Gamma_1\Gamma_2}_{n_1n_2}} : & \!\!\!\! & {\Gamma_1 \subset \mathcal{S}_{n_1}} &\!\!\!{\land}\!\!\!& {\Gamma_2 \subset \mathcal{S}_{n_2}} &&{\!\!\imp\!\!}& {\Gamma_1 \cup \Gamma_2 \subset \mathcal{S}_{\mathit{max}(n_1,n_2)}}\\
\hjoin^{\Gamma_1\Gamma_2}_{nn} & (\!\!\!\! &f_1&\!\!\!,\!\!\!& f_2 &\!\!\!\!)&\!\!\defeq\!\!& \join^{\Gamma_1\Gamma_2}_{n} (f_1, f_2) &\\
\hjoin^{\Gamma_1\Gamma_2}_{n_1'+1>n_2} & (\!\!\!\! &\IS(f_1)&\!\!\!,\!\!\!& f_2 &\!\!\!\!)&\!\!\defeq\!\!& \IS (\hjoin^{\Gamma_1\Gamma_2}_{n_1'n_2} (f_1, f_2)) & \\
\hjoin^{(\Gamma_1A_1)\Gamma_2}_{n_1'+1>n_2} & (\!\!\!\! &\Iforall(f_1)&\!\!\!,\!\!\!& f_2 &\!\!\!\!)&\!\!\defeq\!\!& \Iforall (\hjoin^{\Gamma_1\Gamma_2}_{n_1'n_2} (f_1, f_2)) & \\
\hjoin^{(\Gamma_1A_1)\Gamma_2}_{n_1'+1>n_2} & (\!\!\!\! &\Iimp(f_1,k_1)&\!\!\!,\!\!\!& f_2 &\!\!\!\!)&\!\!\defeq\!\!& \Iimp (\hjoin^{\Gamma_1\Gamma_2}_{n_1'n_2} (f_1, f_2), k_1) & \\
\hjoin^{\Gamma_1\Gamma_2}_{n_1<n_2'+1} & (\!\!\!\! &f_1&\!\!\!,\!\!\!& \IS(f_2) &\!\!\!\!)&\!\!\defeq\!\!& \IS (\hjoin^{\Gamma_1\Gamma_2}_{n_1n_2'} (f_1, f_2)) & \\
\hjoin^{\Gamma_1(\Gamma_2A_2)}_{n_1<n_2'+1} & (\!\!\!\! &f_1&\!\!\!,\!\!\!&\mathtt{I_\forall}(f_2) &\!\!\!\!)&\!\!\defeq\!\!& \mathtt{I_\forall} (\hjoin^{\Gamma_1\Gamma_2}_{n_1n_2'} (f_1, f_2)) & \\
\hjoin^{\Gamma_1(\Gamma_2A_2)}_{n_1<n_2'+1} & (\!\!\!\! &f_1&\!\!\!,\!\!\!&\mathtt{I_\imp}(f_2,k_2) &\!\!\!\!)&\!\!\defeq\!\!& \mathtt{I_\imp} (\hjoin^{\Gamma_1\Gamma_2}_{n_1n_2'} (f_1, f_2), k_2) & \\
\end{array}
$$

In particular, it has to be noticed that the merge possibly does
arbitrary choices: when the same implicative formula $A$ occurs in both contexts, only one of the two
proofs telling how to reduce $\Gamma, A \vdash \dbot$ to
$\mathcal{T}_0, \dneg A_0 \vdash \dbot$ (third clause of $\join_n$) is (arbitrarily) kept.

Another combinatoric lemma is that the merge of contexts indeed
produces a bigger context. To state the lemma, we already need to
define the inclusion of contexts $\Gamma \subset \Gamma'$. This can be
done inductively by the following clauses:
$$
\seqr{L^{\subset}_0}{}{\epsilon \subset \epsilon}
\qquad
\seqr{L^{\subset}_N}{\Gamma \subset \Gamma'}{\Gamma \subset \Gamma', A}
\qquad
\seqr{L^{\subset}_S}{\Gamma \subset \Gamma'}{\Gamma, A \subset \Gamma', A}
$$

Two straightforward lemmas are that $\Gamma_1 \subset \Gamma_1,
\Gamma_2$ and $\Gamma_2 \subset \Gamma_1, \Gamma_2$ for $\Gamma_1$ and $\Gamma_2$
included in $\mathcal{T}_0$. The proofs are by induction on $\Gamma_1$ and $\Gamma_2$ where $i$ is either 1 or 2:
$$
\begin{array}{lcll}
\inclctx_i^{\Gamma_1,\Gamma_2} & : & {\Gamma_i \subset \Gamma_1, \Gamma_2}\\
\inclctx_i^{\epsilon,\epsilon} & \defeq & L_0\\
\inclctx_1^{(\Gamma_1,A),\epsilon} & \defeq & L_S(\inclctx_1^{\Gamma_1\epsilon})\\
\inclctx_2^{(\Gamma_1,A),\epsilon} & \defeq & L_N(\inclctx_2^{\Gamma_1\epsilon})\\
\inclctx_1^{\Gamma_1,(\Gamma_2,A)} & \defeq & L_N(\inclctx_1^{\Gamma_1\Gamma_2})\\
\inclctx_2^{\Gamma_1,(\Gamma_2,A)} & \defeq & L_S(\inclctx_2^{\Gamma_1\Gamma_2})\\
\end{array}
$$
This allows to prove the following lemma where $i$ is either $1$ or $2$:
$$
\begin{array}{lcll}
\inclcupctx_i^{\Gamma_1,\Gamma_2} & : & {\Gamma_i \subset \Gamma_1 \cup \Gamma_2}\\
\inclcupctx_i^{(\Gamma_1,\dneg A_0),(\Gamma_2,\dneg A_0)} & \defeq & L_S(\inclctx_i^{\Gamma_1,\Gamma_2})\\
\inclcupctx_1^{(\Gamma_1,\dneg A_0),(\Gamma_2,A)} & \defeq & L_N(\inclcupctx_1^{(\Gamma_1,\dneg A_0),\Gamma_2})\\
\inclcupctx_2^{(\Gamma_1,\dneg A_0),(\Gamma_2,A)} & \defeq & L_S(\inclcupctx_2^{(\Gamma_1,\dneg A_0),\Gamma_2})\\
\inclcupctx_1^{(\Gamma_1,A),(\Gamma_2,\dneg A_0)} & \defeq & L_S(\inclcupctx_1^{\Gamma_1,(\Gamma_2,\dneg A_0)})\\
\inclcupctx_2^{(\Gamma_1,A),(\Gamma_2,\dneg A_0)} & \defeq & L_N(\inclcupctx_2^{\Gamma_1,(\Gamma_2,\dneg A_0)})\\
\inclcupctx_i^{(\Gamma_1,A),(\Gamma_2,A)} & \defeq & L_S(\inclcupctx_i^{\Gamma_1,\Gamma_2})\\
\inclcupctx_1^{(\Gamma_1,A),(\Gamma_2,B)} & \defeq & L_N(\inclcupctx_1^{(\Gamma_1,A),\Gamma_2}) & \mbox{if $\godel{A} < \godel{B}$} \\
\inclcupctx_2^{(\Gamma_1,A),(\Gamma_2,B)} & \defeq & L_S(\inclcupctx_2^{(\Gamma_1,A),\Gamma_2}) & \mbox{if $\godel{A} < \godel{B}$} \\
\inclcupctx_1^{(\Gamma_1,A),(\Gamma_2,B)} & \defeq & L_S(\inclcupctx_1^{\Gamma_1,(\Gamma_2,B)}) & \mbox{if $\godel{A} > \godel{B}$} \\
\inclcupctx_2^{(\Gamma_1,A),(\Gamma_2,B)} & \defeq & L_N(\inclcupctx_2^{\Gamma_1,(\Gamma_2,B)}) & \mbox{if $\godel{A} > \godel{B}$} \\
\end{array}
$$

Our object logic is defined by the rules on Figure~\ref{fig:rules}.
Note that we shall use non standard derived rules. For instance, we
shall not use the rule $\abs^\imp$ and $\abs^\forall$ but instead the
derived rules $\pi^{\dimp}_1$, $\pi^{\dimp}_2$ and $\drinker_y$.

Thanks to the previous lemma, we are able to translate proofs of $A_1
\in \mathcal{S}_\omega$ and $A_2 \in \mathcal{S}_\omega$ living in
possibly two different contexts to eventually live in the union
of the two contexts:
$$
\begin{array}{llclccllll}
\merge : & \!\!\!\!A_1 \in \mathcal{S}_{\omega} &\!\!\!{\land}\!\!\!& {A_2 \in \mathcal{S}_{\omega}} &{\!\!\!\imp\!\!\!}& \exists n\exists \Gamma (\Gamma \subset \mathcal{S}_n \land \Gamma \vdash A_1 \land \Gamma \vdash A_2)\\
\merge\!\!\!\! &\!\!\!\! (n_1,\Gamma_1,f_1,p_1)&\!\!\!\!\!\!\!&(n_2,\Gamma_2,f_2,p_2) &\!\!\!\defeq\!\!\!& \left(\begin{array}{l}\mathit{max}(n_1,n_2),(\Gamma_1 \cup \Gamma_2),\hjoin^{}_{}(f_1,f_2),\\\dweak(\inclcupctx_1^{\Gamma_1,\Gamma_2},p_1),\\\dweak(\inclcupctx_2^{\Gamma_1,\Gamma_2},p_2)\end{array}\right)\\
\end{array}
$$
where $\dweak$ is an admissible rule of the object logic.

Thanks to the ability to ensure distinct proofs to live in the same
context, we can reformulate the relevant rules of the object logic as
rules over $\mathcal{S}_\omega$, as well as provide proofs of specific
formulae:
$$
\begin{array}{llclcl}
\APPIMP: & \!\!A \dimp B \in \mathcal{S}_\omega &\!\!\!\!{\land}\!\!\!\!& A \in \mathcal{S}_\omega &{\!\!\imp\!\!}& B \in \mathcal{S}_\omega\\
\APPIMP & \!\!q &\!\!\!\!& q' &\!\!\defeq\!\!& \destlinetopin{\merge\,(q,q')}{(n',\Gamma',f',p'p'')}{(n',\Gamma',f',\appimp(p',p''))}\\
\\
\APPF: & \!\!\dforall x\,A(x) \in \mathcal{S}_\omega &\!\!\!\!{\land}\!\!\!\!& \Term &{\!\!\imp\!\!}& A(t) \in \mathcal{S}_\omega\\
\APPF & \!\!(n,\Gamma,f,p) &\!\!\!\!& t &\!\!\defeq\!\!& (n,\Gamma,f,\appf(p,t))\\
\\
\DNLARGE: & \!\!\dneg\dneg A \in \mathcal{S}_\omega &\!\!\!\!&&{\!\!\imp\!\!}& A \in \mathcal{S}_\omega\\
\DNLARGE & \!\!(n,\Gamma,f,p) &\!\!\!\!&&\!\!\defeq\!\!& (n,\Gamma,f,\mydn(p))\\
\\
\multicolumn{4}{l}{\AXDRINKER_{\dforall x\,A(x)}} &\!\!:\!\!& A(x_n)\dimp\dforall x\,A(x) \in \mathcal{S}_\omega\\
\multicolumn{4}{l}{\AXDRINKER_{\dforall x\,A(x)}} &\!\!\defeq\!\!& (2n+1,(\dneg A_0,A(x_n)\dimp\dforall x\,A(x)),\mathtt{I_\forall}(\inj_{2n}),\ax_0)\\
\multicolumn{4}{l}{} &&\hspace{3.5cm} \mbox{{\it where} $2n = \godel{\dforall x\, A(x)}$}\\
\\
\multicolumn{4}{l}{\AXINIT} &\!\!:\!\!& \dneg A_0 \in \mathcal{S}_\omega\\
\multicolumn{4}{l}{\AXINIT} &\!\!\defeq\!\!& (0,\dneg A_0,J_{\mathit{base}},\ax_0)\\
\end{array}
$$
Similarly, we can formulate rules on provability in $\mathcal{T}_0$:
$$
\begin{array}{llclcl}
\DNABS: & \mathcal{T}_0, \dneg A_0 \vdash \dbot \qquad &{\imp}& \mathcal{T}_0 \vdash A_0\\
\DNABS & (\Gamma,g,p) &\!\!\defeq\!\!& (\Gamma,g,\mydn(\absimp(p)))\\
\\
\LIFTBOT: & \mathcal{T}_0, \dneg A_0 \vdash \dbot &{\imp}& \dbot \in \mathcal{S}_\omega\hspace{5.7cm}\\
\LIFTBOT & (\Gamma,g,p) &\!\!\defeq\!\!& (0,(\Gamma,\dneg A_0),\IO(g),p)\\
\end{array}
$$

\begin{figure}
\begin{center}
{\em Primitive rules}
\end{center}

$$\seqr{\ax_i}{|\Gamma'|=i}{\Gamma, A, \Gamma' \vdash A}
\qquad
\seqr{\mydn}{\Gamma \vdash \dneg \dneg A}{\Gamma \vdash A}
$$
$$
\seqr{\appimp}{\Gamma \vdash A \dimp B \qquad \Gamma' \vdash A}{\Gamma \cup \Gamma' \vdash B}
\qquad
\seqr{\appf_t}{\Gamma \vdash \dforall x\,A(x)}{\Gamma \vdash A(t)}
$$
$$
\seqr{\abs^\imp}{\Gamma, A \vdash B}{\Gamma \vdash A \dimp B}
\qquad
\seqr{\abs^\forall}{\Gamma \vdash A(y) \qquad \mbox{$y$ not in $\dforall x\, A(x), \Gamma$}}{\Gamma \vdash \dforall x\,A(x)}
$$
\smallskip

\begin{center}
{\em Admissible rules}
\end{center}

$$
\seqr{\drinker_y}{\Gamma, A(y)\dimp \dforall x\,A(x) \vdash \dbot\qquad \mbox{$y$ not in $\dforall x\, A(x), \Gamma$}}{\Gamma \vdash \dbot}
$$
$$\seqr{\pi^{\dimp}_1}{\Gamma, A\dimp B \vdash \dbot}{\Gamma \vdash A}
\qquad
\seqr{\pi^{\dimp}_2}{\Gamma, A\dimp B \vdash \dbot}{\Gamma \vdash \dneg B}
\qquad
\seqr{\efq}{\Gamma \vdash \dbot}
            {\Gamma \vdash A}
\qquad
\seqr{\dweak}{\Gamma\subset \Gamma' \qquad \Gamma \vdash A}{\Gamma' \vdash A}
$$
\caption{Inference rules characterising classical first-order predicate calculus}
\label{fig:rules}
\end{figure}

When $\godel{A \dimp B} = 2n+1$, we can also derive the following properties:
$$
\begin{array}{llclcl}
\AXIMP_{A\dimp B}: & \!\!((\mathcal{S}_{2n+1}, A\dimp B \vdash \dbot) \imp (\mathcal{T}_0, \dneg A_0 \vdash \dbot)) &{\!\!\imp\!\!}& A \dimp B \in \mathcal{S}_\omega\\
\AXIMP_{A \dimp B} & \!\!k & \defeq & \left(\begin{array}{l}2n+2,(\dneg A_0, A\dimp B),\\\Iimp(\inj_{2n+1},k),\ax_0\end{array}\right)\\
\\
\PROJLEFT{A \dimp B}: & \!\!\mathcal{S}_{2n+1}, A \dimp B \vdash \dbot &{\!\!\imp\!\!}& A \in \mathcal{S}_\omega\\
\PROJLEFT{A \dimp B} & \!\!(\Gamma,f,p) &\!\!\defeq\!\!& (2n+1,\Gamma,f,\pi^{\dimp}_1\,p)\\
\\
\PROJRIGHT{A \dimp B}: & \!\!\mathcal{S}_{2n+1}, A \dimp B \vdash \dbot &{\!\!\imp\!\!}& \dneg B \in \mathcal{S}_\omega\\
\PROJRIGHT{A \dimp B} & \!\!(\Gamma,f,p) &\!\!\defeq\!\!& (2n+1,\Gamma,f,\pi^{\dimp}_2\,p)\\
\end{array}
$$

We are now ready to present the main computational piece of the
completeness proof and we shall use for that notations reminiscent
from semantic normalisation~\cite{Coquand02}, or type-directed partial
evaluation~\cite{Danvy96}, as considered when proving completeness of
intuitionistic logic with respect to models such a Kripke or Beth
models.

We have to prove $\ptruth{A}{\sigma}{\mathcal{M}_0} \myiff A[\sigma] \in \mathcal{S}_{\omega}$,
which means proving $\ptruth{A}{\sigma}{\mathcal{M}_0} \imp A[\sigma]
\in \mathcal{S}_{\omega}$ and $A[\sigma] \in \mathcal{S}_{\omega} \imp
\ptruth{A}{\sigma}{\mathcal{M}_0}$. As in semantic normalisation (see
Section~\ref{sec:semantic-normalisation}), we shall call \emph{reification}
and write $\downarrow_\sigma^A$ the proof mapping a semantic formula
(i.e. $\ptruth{A}{\sigma}{\mathcal{M}_0}$) to a syntactic formula, i.e.  $A[\sigma]
\in \mathcal{S}_{\omega}$. We shall call \emph{reflection} and write
$\uparrow_\sigma^A$ for the way up going from the syntactic view to the
semantic view.

$$
\begin{array}{llcl}
{\downarrow_\sigma^A~:} & \ptruth{A}{\sigma}{\mathcal{M}_0}  & {\imp} & {A[\sigma]\in \mathcal{S}_\omega}\\
\downarrow_\sigma^{P(\vec{t})} & m & \!\defeq\! & m\\
\downarrow_\sigma^{\dbot} & m & \!\defeq\! & \LIFTBOT(m)\\
\downarrow_\sigma^{A \dimp B} & m & \!\defeq\! & \AXIMP_{A[\sigma] \dimp B[\sigma]}(\kontimp{A\dimp B}{\sigma}(m))\\
\downarrow_\sigma^{\dforall x\, A} & m & \!\defeq\! &
  \APPIMP(\AXDRINKER_{(\forall x\,A)[\sigma]},\downarrow_{\sigma,x\leftarrow x_n}^{A}\! (m\, x_n))~ \mbox{{\it where} $2n = \godel{(\dforall x\, A)[\sigma]}$}\\
\\
{\uparrow_\sigma^A~:} & {A[\sigma] \in \mathcal{S}_\omega} & {\imp} & \ptruth{A}{\sigma}{\mathcal{M}_0}\\
\uparrow_\sigma^{P(\vec{t})} & q & \!\defeq\! & q \\
\uparrow_\sigma^{\dbot} & q & \!\defeq\! & \squeeze\,q\\
\uparrow_\sigma^{A \dimp B} & q & \!\defeq\! & m \mapsto \uparrow_\sigma^B\!(\APPIMP(q,\downarrow_\sigma^A\! m))\\
\uparrow_\sigma^{\dforall x\, A} & q & \!\defeq\! & t \mapsto \uparrow_{\sigma,x\leftarrow t}^{A}\!(\APPF(q,t))\\
\end{array}
$$
where, for $m$ proving $\ptruth{A \dimp B}{\sigma}{\mathcal{M}_0}$,
the relative consistency proof $\kontimp{A \dimp B}{\sigma}(m)$ is defined by:
$$
\!\begin{array}{lcl}
\kontimp{A\dimp B}{\sigma}(m) &\!\!\!\!:\!\!\!\!& (\mathcal{S}_{\godel{A[\sigma]\dimp B[\sigma]}}, A[\sigma] \dimp B[\sigma] \vdash \dbot) \imp (\mathcal{T}_0, \dneg A_0 \vdash \dbot)\\
\kontimp{A\dimp B}{\sigma}(m) &\!\!\!\!\defeq\!\!\!\!& r \!\mapsto\!
\squeeze\,(\APPIMP(\PROJRIGHT{A[\sigma] \dimp B[\sigma]}(r),\downarrow_\sigma^B\! (m\, (\uparrow_\sigma^A\! (\PROJLEFT{A[\sigma] \dimp B[\sigma]}(r))))))\\
\end{array}
$$

We still have to prove that the model is classical, which we do by
lifting the double-negation elimination rule to the semantics:
$$
\begin{array}{lclclcll}
{\class_0} & {:} & \forall A\,((\ptruth{\dneg\dneg A}{\substid}{\mathcal{M}_0}) \imp (\ptruth{A}{\substid}{\mathcal{M}_0}))\\
\class_0 & \defeq & A \mapsto m \mapsto \uparrow_{\substid}^A (\DNLARGE\,(\downarrow_{\substid}^{\dneg\dneg A}m))\\
\end{array}\qquad\qquad
$$

It remains also to show that every formula of $\mathcal{T}_0$ is true
in $\mathcal{M}_0$ and this is obtained by the axiom rule:
$$
\begin{array}{lllcll}
\init_0 & {:} &  \forall B\in \mathcal{T}_0\, \ptruth{B}{\substid}{\mathcal{M}_0}\\
\init_0 & \defeq & B \mapsto h \mapsto \uparrow_{\substid}^B (0,(B,\dneg A_0),\IO(J_{cons}(J_{\mathit{base}},h)),\ax_1)\\
\end{array}\qquad\qquad
$$

Finally, we get the completeness result stated as S1' by:
$$
\begin{array}{lcl}
{\completeness} & {:} & {\forall \mathcal{M}\forall \sigma\, (\ptruth{\classic}{\sigma}{\mathcal{M}} \imp \ptruth{\mathcal{T}_0}{\sigma}{\mathcal{M}}\imp\ptruth{A_0}{\sigma}{\mathcal{M}})} \imp {\mathcal{T}_0 \vdash A_0}\\
\completeness & \defeq & \psi \mapsto \DNABS\,(\squeeze\,(\APPIMP(\AXINIT,\downarrow_\substid^{A_0} (\psi\,\mathcal{M}_0\;\substid\;{\class_0}\;\init_0))))\\
\end{array}
$$
Notice that the final result is a
triple $(\Gamma,g,p)$ such that $p$ is a proof of $\Gamma \vdash A_0$
and $g$ is a proof of $\Gamma \subset \mathcal{T}_0$.

\subsection{The computational content on examples}

To illustrate the behaviour of the completeness proofs, we consider
two examples. We use notations of $\lambda$-calculus to
represent proofs in the meta-logic and constructors from
Figure~\ref{fig:rules} for proofs in the object logic.

We place ourselves in the empty theory and
consider the formula $A_0 \defeq X \dimp Y \dimp X$ with $X$ and $Y$
propositional atoms.

The expansion of $\pvDash A_0$ is $\forall \mathcal{M}\;\forall \sigma\;
(\ptruth{\classic}{\mathcal{M}}{\sigma} \imp \ptruth{X}{}{\mathcal{M}} \imp
\ptruth{Y}{}{\mathcal{M}} \imp \ptruth{X}{}{\mathcal{M}})$. It has a canonical
proof, which, as a $\lambda$-term, is the closure of the so-called combinator
$K$ over the symbols it depends on:
$$m \defeq (\mathcal{D},\mathcal{F},\mathcal{P},B) \mapsto \sigma \mapsto
c \mapsto (x:\mathcal{P}(X)) \mapsto (y:\mathcal{P}(Y)) \mapsto x$$

Applying completeness means instantiating the model by the syntactic
model and the substitution by the empty substitution so as to obtain from $m$ the proof
$$m_0 \defeq (x:X\in\mathcal{S}_\omega) \mapsto (y:Y\in \mathcal{S}_\omega) \mapsto x$$
Our object proof is then the result of evaluating
$$
\DNABS\,(\squeeze\,(\APPIMP(\AXINIT,\downarrow^{A_0} m_0)))
$$
which proceeds as follows:
$$
\DNABS\,(\squeeze\,(\APPIMP(\AXINIT,\AXIMP_{A_0}(\kontimp{A_0}{}(m_0)))))\\
$$
where $\kontimp{A_0}{}(m_0)(r)$
reduces to $\squeeze\,(\APPIMP(\PROJRIGHT{A_0}(r),\downarrow^{Y \dimp X} (m_0\, (\uparrow^X (\PROJLEFT{A_0}(r))))))$
for $r$ proving $\mathcal{S}_{\godel{A_0}}, A_0 \vdash \dbot$, that is to
$$
\squeeze\,(\APPIMP(\PROJRIGHT{A_0}(r),\AXIMP_{Y \dimp X}(\kontimp{Y \dimp X}{}(m_0\,(\PROJLEFT{A_0}(r))))))
$$
In there, $\kontimp{Y \dimp X}{}(m_0\,(\PROJLEFT{A_0}(r)))(r')$, for $r'$ proving $\mathcal{S}_{\godel{Y \dimp X}}, Y
\dimp X \vdash \dbot$, reduces in turn to
$$
\squeeze\,(\APPIMP(\PROJRIGHT{Y \dimp X}(r'),m_0\, (\PROJLEFT{A_0}(r))\, (\PROJLEFT{Y \dimp X}(r'))))
$$
Evaluating $\APPIMP(\AXINIT,\AXIMP_{A_0}(\kontimp{A_0}{}(m_0)))$ gives a tuple
$$(\godel{A_0},(\dneg A_0,A_0),I_\imp(\inj_{\godel{A_0}},\kontimp{A_0}{}(m_0)),p_0)$$
where $p_0 \defeq \appimp(\ax_1,\ax_0)$ is a proof of $\dneg A_0, A_0
\vdash \dbot$ obtained by application of the two axiom rules proving
$\dneg A_0 \vdash \dneg A_0$ and $\dneg A_0, A_0 \vdash A_0$.

Evaluating the outermost $\squeeze$ triggers the application of the
continuation $\kontimp{A_0}{}(m_0)$ to $r_0 \defeq (\dneg
A_0,\inj_{\godel{A_0}},p_0)$, meaning that the whole object proof becomes
$$\DNABS\,(\squeeze\,(\APPIMP(\PROJRIGHT{A_0}(r_0),\AXIMP_{Y \dimp
  X}(\kontimp{Y \dimp X}{}(m_0\,(\PROJLEFT{A_0}(r_0)))))))
$$
Evaluating $\APPIMP(\PROJRIGHT{A_0}(r_0),\AXIMP_{Y \dimp  X}(\kontimp{Y \dimp X}{}(m_0\,(\PROJLEFT{A_0}(r_0)))))$ gives a tuple
$$(\godel{Y \dimp X},(\dneg A_0,Y \dimp X),I_\imp(\inj_{\godel{Y \dimp X}},\kontimp{Y \dimp X}{}(m_0\,(\PROJLEFT{A_0}(r_0)))),p_1)$$
where $p_1 \defeq \appimp(\pi^{\dimp}_2\,(p_0),\ax_0)$ is a proof of $\dneg A_0, Y \dimp X
\vdash \dbot$.

Evaluating the new outermost $\squeeze$ triggers in turn the application of the
continuation $\kontimp{Y \dimp X}{}(m_0\,(\PROJLEFT{A_0}(r_0)))$ to
$r_1 \defeq (\dneg A_0,\inj_{\godel{Y \dimp X}},p_1)$ and this results in
\begin{equation}
\label{eqn}
\DNABS\,(\squeeze\,(\APPIMP(\PROJRIGHT{Y \dimp X}(r_1),m_0\, (\PROJLEFT{A_0}(r_0))\, (\PROJLEFT{Y \dimp X}(r_1)))))
\end{equation}
that is, taking into account the definition of $m_0$
$$\DNABS\,(\squeeze\,(\APPIMP(\PROJRIGHT{Y \dimp X}(r_1),\PROJLEFT{A_0}(r_0))))
$$
No continuations are produced by $\APPIMP(\PROJRIGHT{Y \dimp
  X}(r_1),\PROJLEFT{A_0}(r_0))$ so the only role of the last
$\squeeze$ is to peel the $\IS$ leading to a proof $r_2 \defeq
(\epsilon,J_{\mathit{base}},\mydn(p_2))$ where $p_2 \defeq
\appimp(\pi^{\dimp}_2\,(p_1),\pi^{\dimp}_1\,(p_0))$
combines a proof of $\dneg A_0 \vdash \dneg X$ with a proof of $\dneg
A_0 \vdash X$ to get a proof of $\dneg A_0 \vdash \dbot$.

To summarise, the object proof produced is:

$$
\infer[\mydn]
{\vdash A_0}
{\infer[\absimp]
 {\vdash \dneg \dneg A_0}
 {\infer[\appimp]
  {\dneg A_0 \vdash \dbot}
  {\infer[\pi^{\imp}_2]
   {\dneg A_0 \vdash \dneg X}
   {\infer[\appimp]
     {\dneg A_0, Y \dimp X \vdash \dbot}
     {\infer[\pi^{\imp}_2]
       {\dneg A_0 \vdash \dneg (Y \dimp X)}
       {\infer*[p_0]{\dneg A_0, A_0 \vdash \dbot}{}}
       &
       \infer[\ax]
        {\dneg A_0, Y\dimp X \vdash Y \dimp X}
        {}
     }
   }
   &
   \infer[\pi^{\imp}_1]
     {\dneg A_0 \vdash X}
     {\infer*[p_0]{\dneg A_0, A_0 \vdash \dbot}{}}
  }
 }
}
$$

where $p_0$ is:
$$
\infer[\appimp]
       {\dneg A_0, A_0 \vdash \dbot}
       {\infer[\ax]
         {\dneg A_0 \vdash \dneg A_0}
         {}
        &
        \infer[\ax]
         {A_0 \vdash A_0}
         {}
       }
$$

As a matter of comparison, for the canonical proof of the validity of
$A_0' \defeq X \dimp Y \dimp Y$, everything up to step (\ref{eqn})
above is the same modulo the change of $A_0$ into $A_0'$ and of $Y
\dimp X$ into $Y \dimp Y$. After step (\ref{eqn}), one obtains
$$\DNABS\,(\squeeze\,(\APPIMP(\PROJRIGHT{Y \dimp Y}(r_1'),\PROJLEFT{A_0'}(r_1'))))
$$
where
$$
\begin{array}{lll}
r_1' &\defeq& (\dneg A_0',\inj_{\godel{Y \dimp Y}},p_1')\\
p_1' &\defeq& \appimp(\pi^{\dimp}_2\,(p_0'),\ax_0)\\
p_0' &\defeq& \appimp(\ax_1,\ax_0)\\
\end{array}
$$
Finally, this yields $(\epsilon,J_{\mathit{base}},\mydn(p_2'))$ where $p_2' \defeq
\appimp(\pi^{\dimp}_2\,(p_1'),\pi^{\dimp}_1\,(p_1'))$,
that is, graphically:
$$
\infer[\mydn]
{\vdash A_0'}
{\infer[\absimp]
 {\vdash \dneg \dneg A_0'}
 {\infer[\appimp]
  {\dneg A_0' \vdash \dbot}
   {\infer[\pi^{\imp}_2]
     {\dneg A_0' \vdash \dneg Y}
     {\infer*[p_1']{\dneg A_0', Y \dimp Y \vdash \dbot}{}}
    &
    \infer[\pi^{\imp}_1]
      {\dneg A_0' \vdash Y}
      {\infer*[p_1']{\dneg A_0', Y \dimp Y \vdash \dbot}{}}
   }
 }
}
$$
where, graphically, $p_1'$ is:
$$
\infer[]
 {\dneg A_0', Y \dimp Y \vdash \dbot}
 {\infer[\pi^{\imp}_2]
   {\dneg A_0' \vdash \dneg (Y \dimp Y)}
   {\infer[\appimp]
     {\dneg A_0', A_0' \vdash \dbot}
     {\infer[\ax]
       {\dneg A_0' \vdash \dneg A_0'}
       {}
       \!&\!
       \infer[\ax]
             {A_0' \vdash A_0'}
             {}
     }
   }
   &
   \infer[\ax]
         {\dneg A_0', Y\dimp Y \vdash Y \dimp Y}
         {}
 }
$$

We can notice in particular that, treating the metalanguage as a
$\lambda$-calculus as we did, the two canonical proofs of validity of
$X \dimp X \dimp X$ would not produce the same object language proofs.

\subsection{Extension to conjunction}
\label{sec:conjunction}

\newcommand{\dpair}{\dot\objectrulestyle{pair}}
\newcommand{\andpileft}{\dot\objectrulestyle{\pi^{\wedge}_1}}
\newcommand{\andpiright}{\dot\objectrulestyle{\pi^{\wedge}_2}}
\newcommand{\andpibothsides}[1]{\dot\objectrulestyle{\pi^{\wedge}_{#1}}}

Henkin's original proof~\cite{Henkin49a} includes only implication, universal
quantification and the false connective. Handling conjunction in our presentation of Henkin's proof is
straightforward. Let us assume the object language being equipped with
the following rules for conjunction:
$$
\seqr{\dpair}{\Gamma \vdash A_1 \quad \Gamma \vdash A_2}{\Gamma \vdash A_1 \dand A_2}
\qquad
\seqr{\andpileft}{\Gamma \vdash A_1 \dand A_2}{\Gamma \vdash A_1}
\qquad
\seqr{\andpiright}{\Gamma \vdash A_1 \dand A_2}{\Gamma \vdash A_2}
$$

Truth for conjunction being defined by
$$\ptruth{A_1 \dand A_2}{\sigma}{\mathcal{M}} ~~\defeq~~ \ptruth{A_1}{\sigma}{\mathcal{M}}
\wedge \ptruth{A_2}{\sigma}{\mathcal{M}}
$$ a case is added for
conjunctive formulae in each direction of the proof of
$\truth{A}{\sigma}{\mathcal{M}_0} \myiff A[\sigma] \in \mathcal{S}_{\omega}$ as
follows:
$$
\begin{array}{llcl}
\downarrow_\sigma^{A_1 \dand A_2} & (m_1,m_2) & \defeq & \ANDPAIR{A_1\land A_2}\,(\downarrow_\sigma^{A_1}\!m_1, \downarrow_\sigma^{A_2}\!m_2)\\
\uparrow_\sigma^{A_1 \dand A_2} & q & \defeq & (\uparrow_\sigma^{A_1}\!(\ANDPROJLEFT{A_1\land A_2}\,q),\uparrow_\sigma^{A_2}\!(\ANDPROJRIGHT{A_1\land A_2}\,q))\\
\end{array}
$$ where the following combinators lift inference rules on $\mathcal{S}_\omega$:
$$
\begin{array}{lll}
\ANDPAIR{A_1 \dand A_2}(q_1,q_2) &\!\!\defeq\!\!& \dest{\merge(q_1,q_2)}{(n,\Gamma,f,p_1,p_2)}{(n,\Gamma,f,\dpair(p_1,p_2))}\\
\ANDPROJBOTHSIDES{i}{A_1 \dand A_2}(n,\Gamma,f,p) &\!\!\defeq\!\!& (n,\Gamma,f,\andpibothsides{i}\,p)\\
\end{array}
$$
In particular, there is no need to consider
conjunctive formulae in the enumeration.

\subsection{Extension to disjunction}
\label{sec:disjunction}

Taking inspiration from works on
normalisation-by-evaluation in the presence of disjunction
(e.g.~\cite{AbelSattler19,AltDybHofSco01}), we give a proof for
disjunction that relies on the
following arithmetical generalised form of Double Negation Shift
introduced\footnote{In the presence of Markov's principle, this
  generalised form of $\DNS$ is equivalent to the usual form. The
  interest of the generalised form is precisely that it can be used in
  situations which would have required Markov's principle without
  requiring Markov's principle explicitly.} in~\cite{Ilik11}:
$$
\begin{array}{ll}
\DNS^{\forall}_{\mathcal{C}} \quad & \forall n\, ((\mathcal{A}(n) \imp \mathcal{C}) \imp \mathcal{C}) \imp (\forall n\, \mathcal{A}(n) \imp \mathcal{C}) \imp \mathcal{C}
\end{array}
$$
for $\mathcal{C}$ a $\Sigma^0_1$-formula and $\mathcal{A}$ arbitrary (usual $\DNS$ is then $\DNS^\forall_\bot$).

The proof however requires a significant change: instead of proving
$\ptruth{A}{\sigma}{\mathcal{M}_0} \Leftrightarrow {A[\sigma]\in
  \mathcal{S}_\omega}$, we prove:
$$
\begin{array}{lll}
\downarrow^A_\sigma &:& \ptruth{A}{\sigma}{\mathcal{M}_0} \imp {A[\sigma]\in \mathcal{S}_\omega}\\
{\uparrow'_\sigma}^A &:& {A[\sigma]\in \mathcal{S}_\omega} \imp \mathsf{KONT}(\ptruth{A}{\sigma}{\mathcal{M}_0})\\
\end{array}
$$
where $\mathsf{KONT}(\mathcal{A})$, a continuation
monad\footnote{Note that
  $\mathsf{KONT}(\ptruth{A}{\sigma}{\mathcal{M}_0})$ is actually the same
  as $\ptruth{\dneg\dneg A}{\sigma}{\mathcal{M}_0}$ but we shall
  use $\mathsf{KONT}$ also on formulae which are not of the form
  $\ptruth{A}{\sigma}{\mathcal{M}_0}$.}, is defined by:
$$
\mathsf{KONT}(\mathcal{A}) ~~~\defeq~~~ (\mathcal{A} \,\imp\, \mathcal{T}_0, \dneg A_0 \vdash \dbot) \,\imp\,  \mathcal{T}_0, \dneg A_0 \vdash \dbot
$$
As such, the proof can also be seen as an adaptation to Tarski
semantics of the proof of completeness with respect to
continuation-passing style models given in~\cite{Ilik11b}.

Let us assume the object language equipped with the following rules for disjunction:
$$
\seqr{\dinj_1}{\Gamma \vdash A_1}
     {\Gamma \vdash A_1 \dor A_2}
\quad
\seqr{\dinj_2}{\Gamma \vdash A_2}
     {\Gamma \vdash A_1 \dor A_2}
\quad
\seqr{\dcaseraw}
     {\Gamma \vdash A_1 \dor A_2 \quad \Gamma \vdash A_1 \dimp B \quad \Gamma \vdash A_2 \dimp B}
     {\Gamma \vdash B}
$$
Let us also consider the following lifting of the inference
rules to provability in $\mathcal{S}_\omega$:
$$
\begin{array}{llclclllll}
\ORINJBOTHSIDES{i}{A_1\dor A_2}: & A_i \in \mathcal{S}_\omega &{\!\!\imp\!\!}& A_1 \dor A_2 \in \mathcal{S}_\omega\\
\ORINJBOTHSIDES{i}{A_1\dor A_2} & (n,\Gamma,f,p) &\!\!\defeq\!\!& (n,\Gamma,f,\dinj_i\,p)\\
\\
\CASE: & \!\!\!\left(\begin{array}{l}A \dor B \in \mathcal{S}_\omega\\ A \dimp C \in \mathcal{S}_\omega\\ B \dimp C \in \mathcal{S}_\omega\end{array}\right)&{\!\!\imp\!\!}& C \in \mathcal{S}_\omega\\
\CASE & (q,q_1,q_2) & \!\!\defeq\!\!& \destlinetopin{\mergethree(q,q_1,q_2)}{(n,\Gamma,f,p,p_1,p_2)}{(n,\Gamma,f,\dcase{p}{p_1}{p_2})}\\
\end{array}
$$
where we needed the following three-part variant $\mergethree$ of $\merge$:
$$
\!\!\begin{array}{llcll}
\mergethree\! :\! &\!\!\! A_1 \in \mathcal{S}_{\omega} \land A_2 \in \mathcal{S}_{\omega} \land A_3 \in \mathcal{S}_{\omega} &\!\!\!\!\imp\!\!\!\!& \exists n\exists \Gamma (\Gamma \subset \mathcal{S}_n \land \Gamma \vdash A_1 \land \Gamma \vdash A_2 \land \Gamma \vdash A_3)\\
\mergethree &\!\!\!(q_1,q_2,(n,\Gamma,f,p_3)) & \!\!\!\!\defeq\!\!\!\! &
\destlinetopin{\merge(q_1,q_2)}{(n',\Gamma',f',p_1,p_2)}
              {\destlinetwo{\merge((n',\Gamma',f',p_2),(n,\Gamma,f,p_3))}{(n'',\Gamma'',f'',p_2',p_3)}{(n'',\Gamma'',f'',\dweak(\inclcupctx^{\Gamma',\Gamma}_1,p_1),p_2',p_3)}}\\
\end{array}
$$

Truth for disjunction being defined by
$$\ptruth{A_1 \dor A_2}{\sigma}{\mathcal{M}} ~~\defeq~~
\ptruth{A_1}{\sigma}{\mathcal{M}} \vee
\ptruth{A_2}{\sigma}{\mathcal{M}}$$ we show the extended reification
below where $\metarulestyle{case}$ does a case analysis on the form
$\inj_1(m)$ or $\inj_2(m)$ of a proof of disjunction in the typed
$\lambda$-calculus which we use to represent our metalanguage.

$$
\begin{array}{llcl}
{\downarrow_\sigma^A~:} & \ptruth{A}{\sigma}{\mathcal{M}_0}  & {\imp} & {A[\sigma]\in \mathcal{S}_\omega}\\
\downarrow_\sigma^{P(\vec{t})} & m & \!\defeq\! & m\\
\downarrow_\sigma^{\dbot} & m & \!\defeq\! & \LIFTBOT(m)\\
\downarrow_\sigma^{A \dimp B} & m & \!\defeq\! & \AXIMP_{A[\sigma] \dimp B[\sigma]}(\kontimpprim{A\dimp B}{\sigma}(m))\\
\downarrow_\sigma^{\dforall x\, A} & m & \!\defeq\! &
  \APPIMP(\AXDRINKER_{(\forall x\,A)[\sigma]},\downarrow_{\sigma,x\leftarrow x_n}^{A}\! (m\, x_n))~ \mbox{{\it where} $2n = \godel{(\dforall x\, A)[\sigma]}$}\\
\downarrow_\sigma^{A_1 \dand A_2} & (m_1,m_2) & \defeq & \ANDPAIR{A_1\land A_2}\,(\downarrow_\sigma^{A_1}\!m_1, \downarrow_\sigma^{A_2}\!m_2)\\
\downarrow^{A_1 \dor A_2}_{\sigma} &\!\! m\!\! & \!\!\!\defeq\!\!\! & \Ccaseline{m}{\inj_1(m_1)}{\ORINJLEFT{A_1\dor A_2}(\downarrow^{A_1}_{\sigma}\!\!m_1)}{\inj_2(m_2)}{\ORINJRIGHT{A_1 \dor A_2}(\downarrow^{A_2}_{\sigma}\!\!m_2)}\\
\end{array}
$$

Note that this extended reification is unchanged except for the replacement
of $\kontimp{A\dimp B}{\sigma}$ by $\kontimpprim{A\dimp B}{\sigma}$
so as to take into account the use of $\mathsf{KONT}$ in $\uparrow'$:
$$
\kontimpprim{A\dimp B}{\sigma}(m)(r) \defeq
(\uparrow_\sigma^{'A} (\PROJLEFT{A \dimp B}(r)))\, (m' \mapsto
\squeeze\,(\APPIMP(\PROJRIGHT{A \dimp B}(r),\downarrow_\sigma^B (m m'))))
$$

Before giving the modified reflection proof, we need to prove a form
of ex falso quodlibet deriving the truth of any formula $A$ from any
inconsistency $\mathcal{T}_0, \dneg A_0 \vdash \dbot$. This is a
standard proof by induction on $A$:
$$
\begin{array}{lllll}
\EFQ{\sigma}^A : & \mathcal{T}_0, \dneg A_0 \vdash \dbot & \imp & \ptruth{A}{\sigma}{\mathcal{M}_0}\\
\EFQ{\sigma}^{P(\vec{t})} & (\Gamma,g,p) & \defeq & (0,(\Gamma,\dneg A_0),\IO(g),\efq\, p)\\
\EFQ{\sigma}^{\dbot} & (\Gamma,g,p) & \defeq & (\Gamma,g,p)\\
\EFQ{\sigma}^{A \dimp B} & (\Gamma,g,p) & \defeq & m \mapsto \EFQ{\sigma}^{B} (\Gamma,g,p)\\
\EFQ{\sigma}^{\dforall x\,A} & (\Gamma,g,p) & \defeq & t \mapsto \EFQ{\sigma,x\leftarrow t}^{A} (\Gamma,g,p)\\
\EFQ{\sigma}^{A_1 \dand A_2} & (\Gamma,g,p) & \defeq & (\EFQ{\sigma}^{A_1} (\Gamma,g,p), \EFQ{\sigma}^{A_2} (\Gamma,g,p))\\
\EFQ{\sigma}^{A_1 \dor A_2} & (\Gamma,g,p) & \defeq & \inj_1\,(\EFQ{\sigma}^{A_1} (\Gamma,g,p))\\
\end{array}
$$
where we may notice in passing that an arbitrary choice is made in the disjunction case.

We need the shift of double negation with respect to implication to
formulae expressing truth. Ex falso quodlibet being obtained on formulae
expressing truth with respect to possibly-exploding models, the proof
is standard:
$$
\begin{array}{lcl}
\DNS^{\imp} &\!\!:\!\!& (\mathcal{A} \imp \mathsf{KONT}(\ptruth{B}{\sigma}{\mathcal{M}_0})) ~\imp~ \mathsf{KONT}(\mathcal{A} \imp \ptruth{B}{\sigma}{\mathcal{M}_0})\\
\DNS^{\imp} &\!\!\defeq\!\! & H \mapsto K \mapsto K\, (m_A \mapsto \EFQ{\sigma}^B\, (H~m_A~(m_B \mapsto K~(m_A' \mapsto m_B))))\\
\end{array}
$$

We also need the shift of double negation with respect to conjunction. The intuitionistic proof is easy:
$$
\begin{array}{lcl}
\DNS^{\land} &\!\!:\!\!& \mathsf{KONT}(\mathcal{A}_1) \land \mathsf{KONT}(\mathcal{A}_2) \imp \mathsf{KONT}(\mathcal{A}_1 \land \mathcal{A}_2)\\
\DNS^{\land} &\!\!\defeq\!\! & (H_1,H_2) \mapsto K \mapsto H_1\, (m_1 \mapsto H_2\, (m_2 \mapsto K~(m_1,m_2)))\\
\end{array}
$$

We are now ready to reformulate reflection, including the case for disjunction:
$$
\begin{array}{llcl}
{\uparrow_\sigma^{'A~}:} & {A[\sigma] \in \mathcal{S}_\omega} & {\rightarrow} & \mathsf{KONT}\,(\ptruth{A}{\sigma}{\mathcal{M}_0})\\
\uparrow_\sigma^{'P(\vec{t})} & q & \!\defeq\! & K \mapsto K\,q \\
\uparrow_\sigma^{'\dbot} & q & \!\defeq\! & K \mapsto \squeeze\,q\\
\uparrow_\sigma^{'A \dimp B} & q & \!\defeq\! & \DNS^{\imp} (m \mapsto \uparrow_\sigma^{'B} (\APPIMP(q,\downarrow_A m)))\\
\uparrow_\sigma^{'\dforall x\, A(x)} & q & \!\defeq\! & \DNS^{\forall}\, (t \mapsto \uparrow_{\sigma,x\leftarrow t}^{'A}(\APPF(q,t)))\\
\uparrow_\sigma^{'A_1 \dand A_2} & q & \defeq & \DNS^{\land}(\uparrow_\sigma^{'A_1}(\ANDPROJLEFT{A_1 \dand A_2}\,q),\uparrow_\sigma^{'A_2}(\ANDPROJRIGHT{A_1 \dand A_2}\,q))\\
\uparrow_\sigma^{'A_1 \dor A_2} & q & \!\!\defeq\!\! & K \mapsto \squeeze\,(\CASE(q,\!\!\!\begin{array}[t]{l}\AXIMP_{\dneg A_1}(\kontorkont{A_1\dor A_2}{\sigma,1}(K)),\\\AXIMP_{\dneg A_2}(\kontorkont{A_1\dor A_2}{\sigma,2}(K))))\end{array}\\
\end{array}
$$
where, for $K$ proving $\ptruth{\dneg (A_1\dor A_2)}{id}{\mathcal{M}_0}$, the continuation $\kontorkont{A_1\dor A_2}{\sigma,i}(K)$ is defined by:
$$
\begin{array}{lclclcll}
\kontorkont{A_1\dor A_2}{\sigma,i}(K) &\!\!\!:\!\!\!& (\mathcal{S}_{\godel{\dneg A_i[\sigma]}}, \dneg A_i[\sigma] \vdash \dbot) ~\imp~ (\mathcal{T}_0,\dneg A_0 \vdash \dbot)\\
\kontorkont{A_1\dor A_2}{\sigma,i}(K) & \defeq & r \mapsto (\uparrow_\sigma^{'A_i}(\PROJLEFT{\dneg A_i[\sigma]}(r)))\,(m \mapsto K\,(\inj_i\,m))\\
\end{array}
$$

The proof of {$\class_0$} follows a different pattern than the one without $\mathsf{KONT}$. For
it and for the proof of {$\init_0$}, we use again $\DNS^{\forall}$ to distribute the
quantification over the axioms of the theory:
$$
\begin{array}{lllcll}
{\class_0} & {:} & \mathsf{KONT}(\forall A\, ((\ptruth{\dneg\dneg A}{\substid}{\mathcal{M}_0}) \imp (\ptruth{A}{id}{\mathcal{M}_0})))\\
\class_0 & \defeq & \DNS^{\forall}\, (A \mapsto \DNS^{\imp} (m \mapsto m))\\
\\
\init_0 & {:} &  \mathsf{KONT}(\forall B\in \mathcal{T}_0\, \ptruth{B}{\substid}{\mathcal{M}_0})\\
\init_0 & \defeq & \DNS^{\forall}\, (B \mapsto \DNS^{\imp} (h \mapsto \uparrow_\substid^{'B}(0,(B,\dneg A_0),\IO(J_{cons}(J_{\mathit{base}},h)),\ax_1)))\\
\end{array}
$$
Finally, the proof of completeness also needs to chain continuations:
$$
\begin{array}{lllcl}
{\completeness} & {:} & {\forall \mathcal{M}\forall \sigma\, (\ptruth{\classic}{\sigma}{\mathcal{M}} \imp \ptruth{\mathcal{T}_0}{\sigma}{\mathcal{M}}\imp\ptruth{A_0}{\sigma}{\mathcal{M}})} \imp {\mathcal{T}_0 \vdash A_0}\\
\completeness & \defeq & \psi \mapsto \DNABS\,({\completeness'}\,\psi)\\
\end{array}
$$
where
$$
\!\begin{array}{lllcl}
{\completeness'} : {\forall \mathcal{M}\forall \sigma\, (\ptruth{\classic}{\sigma}{\mathcal{M}} \imp \ptruth{\mathcal{T}_0}{\sigma}{\mathcal{M}}\imp\ptruth{A_0}{\sigma}{\mathcal{M}})} \imp \dneg A_0, {\mathcal{T}_0 \vdash \dbot}\\
\completeness' \defeq \psi \!\mapsto\! \class_0\, (c \!\mapsto\! \init_0\, (i \!\mapsto\! \squeeze\,(\APPIMP(\AXINIT,\downarrow_\substid^{A_0} (\psi\,\mathcal{M}_0\,\substid\,c\,i)))))\\
\end{array}
$$

An interesting remark is that $\mathcal{T}_0, \dneg A_0 \vdash \dbot$
intuitionistically implies $\mathsf{KONT}(\bot)$, so, by using
$\mathsf{KONT}$, reflection for $\dbot$ does not need any more to
consider possibly-exploding model in order to be intuitionistically
valid. The need for possibly-exploding models (or Markov's principle)
shows up instead in proving $\EFQ{\sigma}^{\dbot} : \mathcal{T}_0,
\dneg A_0 \vdash \dbot \imp \ptruth{\dbot}{\sigma}{\mathcal{M}_0}$
which is in turn required to prove our specific version of
$\DNS^{\imp}$, itself used to prove reflection for the connective
$\dimp$, as well as to prove $\class_0$ and $\init_0$.

\subsection{About the logical strength of completeness in the presence of disjunction}

Kirst~\cite[Fact 7.36]{Kirst23} proved in the context of
intuitionistic propositional epistemic logic that completeness
with respect to (non-exploding) Kripke semantics implies the following
classical principle:
$$
\neg\neg\forall n\,(\neg \mathcal{A}(n) \lor \neg\neg \mathcal{A}(n))
$$
which in turn is equivalent to the following instance of $\DNS^{\forall}_{\bot}$ for disjunctions:
$$
\DDNS \qquad \forall n\,\neg\neg\,(\neg \mathcal{A}(n) \lor \neg \mathcal{B}(n)) \imp \neg\neg\forall n\,(\neg \mathcal{A}(n) \lor \neg \mathcal{B}(n))
$$

In private communication, Kirst even showed that $\DDNS$ is enough,
showing that it is exactly the extra bit of classical reasoning needed
to handle disjunction in a constructive metalanguage (in the case of
intuitionistic propositional epistemic logic).

Kirst's proof\footnote{A variant of his original proof directly
proving $\DDNS$ from the completeness of
bi-intuitionistic logic with respect to its Kripke semantics can be
found in~\cite{ShillitoKirst24}.}
can be adapted to the case of completeness of
classical logic with respect to non-exploding Tarski semantics. We
give here a variant of his proof showing that completeness in the presence of $\dimp$, $\dor$ and $\dbot$ relative to
enumerable theories and a $\Sigma^0_1$ exploding interpretation of $\dbot$ implies the
following generalised\footnote{Like in the case of
$\DNS^{\forall}_{\mathcal{C}}$, $\DDNS_{\mathcal C}$ is useful to reason in the
absence of Markov's principle but is equivalent to $\DDNS$ in the
presence of Markov's principle.}\,\footnote{Using
conjunction, $\DDNS_{\mathcal C}$ can be equivalently stated as
$\forall n\, \neg_{\mathcal{C}} (\mathcal{A}(n) \land \mathcal{B}(n)) \imp \neg_{\mathcal{C}}(\forall \mathcal{D}\,\neg_{\mathcal{D}}\forall n\,(\neg_{\mathcal{D}} \mathcal{A}(n) \lor \neg_{\mathcal{D}} \mathcal{B}(n)))$.}
form of $\DDNS$ for $\mathcal{A}$, $\mathcal{B}$ and
$\mathcal{C}$ $\Sigma_0^1$-formulas and for $\mathcal{D}$ ranging over
$\Sigma_0^1$-formulas too:
$$
\DDNS_{\mathcal C} \qquad \forall n\,\neg_{\mathcal{C}}(\forall \mathcal{D}\,\neg_{\mathcal{D}}\,(\neg_{\mathcal{D}} \mathcal{A}(n) \lor \neg_{\mathcal{D}} \mathcal{B}(n))) \imp \neg_{\mathcal{C}}(\forall \mathcal{D}\,\neg_{\mathcal{D}}\forall n\,(\neg_{\mathcal{D}} \mathcal{A}(n) \lor \neg_{\mathcal{D}} \mathcal{B}(n)))
$$
where $\neg_{\mathcal{F}} \mathcal{E}$ abbreviates $\mathcal{E} \imp \mathcal{F}$.

Indeed, take a language with families of atoms $X_n$ and $Y_n$,
and consider the theory $\mathcal{T}$
made of $\dneg X_n \dor \dneg Y_n$ for each $n$, as well as $X_n$
for each $n$ such that $\mathcal{A}(n)$, and $Y_n$ for each $n$ such that
$\mathcal{B}(n)$. The idea of the proof is that $\forall \mathcal{D}\,
\neg_{\mathcal{D}} \forall n\,(\neg_{\mathcal{D}} \mathcal{A}(n) \lor
\neg_{\mathcal{D}} \mathcal{B}(n))$ implies\footnote{Indeed, assuming
a model of $\mathcal{T}$ interpreting the atoms $X_n$ and $Y_n$ by
some predicates $\mathcal{X}$ and $\mathcal{Y}$ respectively, as well
as $\dbot$ by some $\Sigma^0_1$ proposition $\mathcal{Z}$, then, the truth of
$X_n$ whenever $\mathcal{A}(n)$ implies
$\neg_{\mathcal{Z}}\mathcal{X}(n) \imp
\neg_{\mathcal{Z}}\mathcal{A}(n)$, and, similarly
$\neg_{\mathcal{Z}}\mathcal{Y}(n) \imp
\neg_{\mathcal{Z}}\mathcal{B}(n)$. Moreover, the truth of $\dneg X_n
\dor \dneg Y_n$ for each $n$ implies
$\neg_{\mathcal{Z}}\mathcal{X}(n) \lor
\neg_{\mathcal{Z}}\mathcal{Y}(n)$, thus
$\neg_{\mathcal{Z}}\mathcal{A}(n) \lor
\neg_{\mathcal{Z}}\mathcal{B}(n)$ for each $n$. We then get
$\mathcal{Z}$ by the initial hypothesis $\forall
\mathcal{D}\neg_{\mathcal{D}}\forall n\,(\neg_{\mathcal{D}}
\mathcal{A}(n) \lor \neg_{\mathcal{D}} \mathcal{B}(n))$. Remark that
restricting the interpretation of $\dbot$ to $\Sigma^0_1$ propositions
is actually not a restriction since the syntactic model in Henkin's
proof for enumerable theory $\mathcal{U}$ and conclusion $A$ precisely
interprets $\dbot$ by $\mathcal{U}, \dneg{A} \vdash \dbot$ which is
$\Sigma^0_1$ and because, by composition with soundness, completeness is equivalent to
its restriction to the syntactic model.} $\mathcal{T}
\vDash \dbot$, which implies $\Gamma \vdash \dbot$ for a finite subset
$\Gamma$ of $\mathcal{T}$ by completeness, which implies in turn
$\Gamma \vDash \dbot$ by soundness. Now, by taking the exploding model
$\mathcal{M}_{\mathit{std}}$ interpreting $X_n$ by $\mathcal{A}(n)$,
$Y_n$ by $\mathcal{B}(n)$ and $\dbot$ by $\mathcal{C}$, the final goal
$\mathcal{C}$, that is $\mathcal{M}_{\mathit{std}} \vDash \dbot$, can eventually be
obtained by proving $\mathcal{M}_{\mathit{std}} \vDash \mathcal{T}$, that is,
$\mathcal{M}_{\mathit{std}} \vDash \dneg X_n \dor \dneg Y_n$, that
is $\neg_{\mathcal{C}} \mathcal{A}(n) \lor \neg_{\mathcal{C}}
\mathcal{B}(n)$, for only those finite number of
$\dneg X_n \dor \dneg Y_n$ in~$\Gamma$. The latter can be obtained by first chaining a finite number of
instances of $\forall
n\,\neg_{\mathcal{C}}\neg_{\mathcal{C}}\,(\neg_{\mathcal{C}}
\mathcal{A}(n) \lor \neg_{\mathcal{C}} \mathcal{B}(n))$, each of them
coming by instantiating $\mathcal{D}$ with $\mathcal{C}$ in $\forall
n\,\neg_{\mathcal{C}}\forall
{\mathcal{D}}\,\neg_{\mathcal{D}}\,(\neg_{\mathcal{D}} \mathcal{A}(n)
\lor \neg_{\mathcal{D}} \mathcal{B}(n))$ for the corresponding $n$'s.

The conjunction of facts that $\DDNS_{\mathcal{C}}$ is provable from
completeness with possibly-exploding semantics, that it is equivalent to
$\DDNS_{\bot}$ in the presence of Markov's principle, that Markov's
principle is the principle needed to handle falsity in a non-exploding
semantics and that Kirst showed that nothing more is needed to handle both
falsity and disjunction suggests that $\DDNS_{\mathcal{C}}$ is exactly
the amount of classical reasoning needed to handle disjunction in the
presence of a possibly-exploding interpretation of falsity.

Even though our proof of completeness with respect to
possibly-exploding models in the presence of disjunction uses a priori the full
generality of $\DNS^{\forall}_{\mathcal{C}}$ for $\exists$-free formulae, we conjecture that
it could be obtained using $\DDNS_{\mathcal{C}}$ instead of $\DNS^{\forall}_{\mathcal{C}}$,
that is, we conjecture that either $\DNS^{\forall}_{\mathcal{C}}$ for
$\exists$-free formulae derives from $\DDNS$ or that the completeness
proof can be modified to depend only on $\DDNS_{\mathcal{C}}$.

Two other recent lines of research about the logical strength of completeness
with disjunction are Krivtsov~\cite{Krivtsov15} and Espíndola~\cite{EspindolaPhD,Espindola16}.
Krivtsov showed that with respect to exploding
models\footnote{Krivtsov calls such models \emph{intuitionistic
structures}.}  G\"odel's completeness for recursively enumerable
theories is equivalent to the Weak Fan Theorem ($\WFT$) over a rather weak intuitionistic arithmetic. This is to be compared
to Veldman's proof of completeness with respect to Kripke
semantics~\cite{Veldman76} which already required $\WFT$. Since
$\WFT$ is classically equivalent to Weak Kőnig's Lemma,
this is also to be compared to the equivalence of G\"odel's
completeness with Weak Kőnig's Lemma over a weak arithmetic with recursive comprehension in the classical reverse
mathematics of the subsystems of second-order
arithmetic~\cite{Simpson99}. On his side, reasoning in the language of
topos theory, Espíndola proved the same result in IZF, where $\WFT$
was formulated as the compactness of the Cantor space.

The situation regarding $\WFT$ is rather subtile as there are
different formulations depending on how infinite paths are represented
in a binary tree (see~\cite{BredeHerbelin21}). Let us call
$\WFT_{\mathit{pred}}$ the intuitionistically provable\footnote{See a
proof relying on intuitionistic $ACA_0$ in the Coq standard
library~\cite[WeakFan.v]{Coqstdlib22}} variant of $\WFT$ where
infinite paths are represented using a predicate and
$\WFT_{\mathit{dec}}$ the weakly classical variant\footnote{See
Berger~\cite{Berger09} who identifies the classical part of (a
functional form of) $\WFT$ as a principle called $L_{\mathit{fan}}$.}
where infinite paths are represented using a decidable predicate. Both
Krivtsov and Espíndola's results refer to $\WFT_{\mathit{dec}}$. On
the other side, Krivine's proof and our own proof for the language
without disjunction are intuitionistic and thus do not require
$\WFT_{\mathit{dec}}$. This is consistent with another result of
Espíndola with Forssell~\cite{Espindola16,ForssellEspindola17} showing
that $\WFT_{\mathit{dec}}$ is not needed on top of IZF to prove
completeness with respect to exploding Kripke semantics in the absence
of disjunction. And indeed, in IZF, like in second-order
intuitionistic arithmetic where Krivine and us are reasoning,
$\WFT_{\mathit{pred}}$ is provable (if we had instead reasoned in a weaker
arithmetic, we would certainly have
explicitly needed $\WFT_{\mathit{pred}}$, the same way as Weak Kőnig's
Lemma was needed in~\cite{Simpson99}). Moreover, this is consistent with yet another
result, namely the formal equivalence between $\WFT_{\mathit{pred}}$
and completeness with respect to Scott entailment relations
in~\cite{BredeHerbelin21}. In turn, we conjecture that
$\WFT_{\mathit{dec}}$, thus completeness with respect to exploding
models and all connectives, disjunction included, is equivalent to
$\WFT_{\mathit{pred}}$ together with $\DDNS_{\mathcal{C}}$.

To conclude, note also that like Markov's principle,
$\DNS^{\forall}_{\mathcal{C}}$ and thus $\DDNS_{\mathcal{C}}$ as well
preserve the witness and disjunction properties of intuitionistic
logic, so they are {\em in this sense} intuitionistically valid. Also,
to compute with $\DNS^{\forall}_{\mathcal{C}}$, bar
recursion~\cite{Spector62} or delimited continuations~\cite{Ilik11}
can be used.

\subsection{Extension to existential quantification}

The situation for existential quantification is simpler than for
disjunction and modifying the statement of reflection is not necessary.  The
idea is to consider an enumeration of formulae which takes existential
formulae into account, then to add a clause to the definition of $\Gamma
\subset \mathcal{S}_n$ similar to the one for universal
quantification, using Henkin's axiom $\dexists y\, A(y) \dimp A(x)$
for $x$ taken fresh in the finite set of formulae coming before
$\dexists y\,A(y)$ in the enumeration\footnote{If universal
  quantification is present among the connectives, we can also reuse
  Henkin's axiom $\dneg A(x) \dimp \dforall y\,\dneg A(y)$ up to some
  extra classical reasoning in the object language.}. Reification is
direct, using the witness coming from the proof of truth as a witness
for the proof of derivability. For reflection, the idea is to combine
a proof of $(\exists y\, A)[\sigma] \in \mathcal{S}_\omega$ with the
proof of $(\dexists y\, A)[\sigma] \dimp A[\sigma,y\leftarrow x]$
available at some level $\mathcal{S}_n$ to get a proof of
$A[\sigma,y\leftarrow x] \in \mathcal{S}_\omega$, then a proof of
$\ptruth{A}{\sigma,y\leftarrow x}{\mathcal{M}}$, thus a proof
\footnote{One may wonder if the ability to extend Henkin's proof to existential
quantification without modifying the statement of reflection is compatible with
the ability to encode $A \dor B$ as an existential formula $\exists
b\,((b=0 \dimp A) \dand (b=1 \dimp B))$ in any signature containing
symbols $0$ and $1$. The point is, that to justify this encoding, we
still need to do a case analysis on the value of $b$ to determine which of
the conjunct we may use in $(b=0 \dimp A) \dand (b=1 \dimp B)$.
This typically requires an axiom of the form $b=0 \dor b=1$, which
moves us back to a situation where $\dor$ is actually already part of the
language.}
of
$\ptruth{\exists y\,A}{\sigma}{\mathcal{M}}$.

\bibliographystyle{asl}

\end{document}